\documentclass{amsart}

\usepackage{mathrsfs,amsfonts,amsmath,amsthm,amssymb}
\usepackage{graphics,pstricks}
\usepackage[all]{xy}
\xyoption{curve}
\xyoption{import}
\xyoption{arc}
\xyoption{ps}
\UsePSspecials{dvips}
\usepackage{graphicx}

\title[Pseudodifferential calculus on manifolds with fibred corners]{Pseudodifferential calculus on manifolds \\ with fibred corners :
         \\ the groupoid of phi-calculus}
\author{Laurent Guillaume}

\DeclareMathOperator{\id}{Id}
\DeclareMathOperator{\im}{Im}

\DeclareMathOperator{\coker}{Coker}
\DeclareMathOperator{\ind}{ind}
\DeclareMathOperator{\symb}{symb}

\DeclareMathOperator{\vol}{vol}

\begin{document}

\hyphenation{mani-fold}

{\theoremstyle{definition}\newtheorem{definition}{Definition}[section]
\newtheorem{notation}[definition]{Notation}
\newtheorem{remnot}[definition]{Remarks and notation}
\newtheorem{terminology}[definition]{Terminology}
\newtheorem{remark}[definition]{Remark}
\newtheorem{remarks}[definition]{Remarks}
\newtheorem{example}[definition]{Example}
\newtheorem{examples}[definition]{Examples}}
\newtheorem{proposition}[definition]{Proposition}
\newtheorem{lemma}[definition]{Lemma}
\newtheorem{theorem}[definition]{Theorem}
\newtheorem{corollary}[definition]{Corollary}

\newtheorem*{dem}{Demonstration}

\newenvironment{demo}{\begin{dem}}%
{\flushright$\qed$\end{dem}}

\newcommand{\ci}{C^{\infty}}
\newcommand{\co}{C^{0}}
\newcommand{\A}{\mathscr{A}}
\newcommand{\B}{\mathscr{B}}
\newcommand{\Cat}{\mathscr{C}}
\newcommand{\D}{\mathscr{D}}
\newcommand{\E}{\mathscr{E}}
\newcommand{\F}{\mathcal{F}_M}
\newcommand{\gr}{\mathscr{G}}
\newcommand{\go}{\mathscr{G} ^{(0)}}
\newcommand{\hr}{\mathscr{H}}
\newcommand{\ho}{\mathscr{H} ^{(0)}}
\newcommand{\gd}{\mathscr{G}^{\mathbb{R}^2}}
\newcommand{\gt}{\mathscr{G} ^{T}}
\newcommand{\I}{\mathscr{I}}
\newcommand{\J}{\mathcal{J}}
\newcommand{\Nb}{\mathscr{N}}
\newcommand{\Kom}{\mathscr{K}}
\newcommand{\ops}{\mathscr{O}}
\newcommand{\Pb}{\mathscr{P}}
\newcommand{\sw}{\mathscr{S}}
\newcommand{\Uo}{\mathscr{U}}
\newcommand{\Vo}{\mathscr{V}}
\newcommand{\Rr}{\mathbb{R}}
\newcommand{\Rn}{\mathbb{R}^n}
\newcommand{\Rp}{\mathbb{R}^p}
\newcommand{\Rq}{\mathbb{R}^q}
\newcommand{\Nat}{\mathbb{N}}
\newcommand{\Rnk}{\mathbb{R}^k_+\times\mathbb{R}^{n-k}}
\newcommand{\Z}{\mathbb{Z}}
\newcommand{\src}{\mathscr{S}_{c}}
\newcommand{\cc}{C_{c}^{\infty}}
\newcommand{\cg}{C_{c}^{\infty}(\gr)}
\newcommand{\cgr}{C_{r}^{\ast}(\gr)}
\newcommand{\cgo}{C_{c}^{\infty}(\go)}
\newcommand{\ct}{C_{c}^{\infty}(\gr^T)}
\newcommand{\ckt}{C_{c}^{k}(\gr \times [0,1])}
\newcommand{\ck}{C_{c}^{k}(\gr)}
\newcommand{\ca}{C_{c}^{\infty}(A\gr)}
\newcommand{\Un}{{U}^{(n)}}
\newcommand{\Du}{D_{\mathscr{U}}}
\newcommand{\X}{\underline{X}}
\newcommand{\mX}{\tilde{X}}
\newcommand{\bX}{\partial{X}}
\newcommand{\bM}{\partial{M}}
\newcommand{\grVF}{\gr V_\mathcal{F}}
\newcommand{\Hol}{\mathcal{H}ol(X,\mathcal{F})}
\newcommand{\HolM}{\mathcal{H}ol(M,\mathcal{F}_M)}
\newcommand{\eq}{\mathscr{E}}
\providecommand{\norm}[1]{\lVert#1\rVert}
\providecommand{\abs}[1]{\lvert#1\rvert} 
\newcommand*{\avint}{\mathop{\ooalign{$\int$\cr$-$}}}






\begin{abstract}
This paper is concerned with pseudodifferential calculus on mani\-folds with fibred corners. Following work of Connes, Monthubert, Skandalis and Androulidakis, we associate to every manifold with fibred corners a longitudinally smooth groupoid which algebraic and differential structure is explicitely described.
This groupoid has a natural geometric meaning as a holonomy groupoid of singular foliation, it is a singular leaf space in the sense of Androulidakis and Skandalis.
We then show that the associated compactly supported pseudo\-differential calculus coincides with Mazzeo and Melrose's $\phi$-calculus and we introduce an extended algebra of smoothing operators that is shown to be stable under holomorphic functional calculus.
This result allows the interpretation of  {$\phi$-calculus} as the pseudo\-differential calculus associated with the holonomy groupoid of the singular foliation defined by the manifold with fibred corners. It is a key step to set the index theory of those singular manifolds in the noncommutative geometry framework.

\end{abstract}

\maketitle

\section*{Introduction}


In order to generalize the Atiyah-Patodi-Singer index theorem, Melrose introduced in \cite{Mel} a pseudodifferential calculus on manifolds with boundary and manifolds with corners : the $b$-calculus. Melrose and Mazzeo then analyzed the case of fibred boundaries and introduced $\phi$-calculus to extend the theorem to families of operators \cite{MM}. 
Moreover the study initiated by Thom \cite{Thom} and Mather \cite{Mather} of pseudo-stratified manifolds and the existence of a desingularization process (\cite{Verona},\cite{BHS}) show that any pseudo-stratified manifold is a quotient space of a manifold with fibred corners (\cite{DL},\cite{DLR}).
Those results motivate the need to understand pseudodifferential calculus on mani\-folds with fibred corners.

The aim of this paper is to set the index theory of those singular manifolds in the noncommutative geometry framework.
In particular we want to obtain the analogue of the leaf space of a foliation for those spaces.
Ehresmann, Haefliger \cite{Haef} and Wilkelnkemper \cite{Win} introduced groupoids to model the leaf space of a regular foliation. Pradines and Bigonnet \cite{BP}, Debord \cite{Deb}, Androulidakis and Skandalis \cite{ASk} studied the case of singular foliations.
Among others they have contributed to realizing the idea that groupoids are natural substitutes to singular spaces.

Following fundamental work by Connes \cite{Coinc,Cosf,Concg}, it also appears that groupoids are key objects to understand pseudodifferential calculus and $K$-theory groups which are the receptacles of the index.
Pseudodifferential calculus on Lie groupoids has been defined independently by Monthubert-Pierrot \cite{MP} and by Nistor, Weinstein and Xu \cite{NWX}. Androulidakis and Skandalis proposed a general framework to deal with singular foliations \cite{ASk2}. Some particular groupoid models for manifolds with
corners and conic manifolds were proposed by Monthubert \cite{Mont} and Debord-Lescure-Nistor \cite{DLN}.
Monthubert showed in \cite{Mont} that to every manifold with corners X could be associated a longitudinally smooth groupoid $\Gamma(X)$ whose pseudodifferential calculus $\Psi^{\infty}(\Gamma(X))$ coincides with Melrose's $b$-calculus. This original result gave new proofs for the study of $b$-calculus.

In this article we show that such a construction exists for $\phi$-calculus by associating to every manifold with fibred boundary, then to every manifold with fibred corners a longitudinally smooth groupoid $\Gamma_\phi(X)$. 
We then show that the associated compactly supported pseudodifferential calculus coincides with Mazzeo and Melrose's $\phi$-calculus and we introduce an extended algebra of smoothing operators that is shown to be stable under holomorphic functional calculus.
Finally we show that the groupoid we built has a natural geometric meaning as a holonomy groupoid of singular foliation, it is an explicit example of a singular leaf space in the sense of Androulidakis and Skandalis \cite{ASk}.
This result allows the interpretation of {$\phi$-calculus} as the pseudodifferential calculus associated with the holonomy groupoid of the singular foliation defined by the manifold with fibred corners.
The reward of this conceptual approach is a simplified exposition of $\phi$-calculus definition and properties, as well as a new geometric interpretation as a holonomy groupoid. Those results are based on the PhD thesis of the author (see \cite{Gui}).

\smallskip 

The paper is organized in 5 sections.
In the first section we recall classical material on groupoids, their associated pseudodifferential calculus and normal cone deformations. These are crucial notions to set the index theory of manifolds with fibred corners in the noncommutative geometry framework. 

In the second section we introduce some definitions related to manifolds with corners and manifolds with embedded fibred corners : following Monthubert the latter are embedded in smooth manifolds endowed with a transverse family of fibred submanifolds of codimension 1 : we call such a structure a \textsl{fibred decoupage}.

In the third section we define the groupoid associated to a fibred decoupage in two steps:
\begin{itemize}
\item we construct a ``puff groupoid" $\grVF$ for any foliated codimension 1 submanifold $(V,\mathcal{F})$ of the manifold M. This groupoid is described as a gluing between the holonomy groupoid of $M\setminus V$ and a deformation groupoid $\D_\varphi = D_\varphi\rtimes\Rr_+^*$:
$$\gr V_\mathcal{F} = \D_\varphi \bigcup_\Psi \mathcal{H}ol(M\setminus V) .$$
$D_\varphi$ is the normal cone deformation of the groupoid immersion $\mathscr{H}ol(\mathcal{F})\rightarrow V\times V$.

\item the groupoid of a fibred decoupage is given by the fibered product of the puff groupoids $\gr V_{\mathcal{F}i}$ for each face.
\end{itemize}

\noindent We can then define the groupoid $\Gamma_\phi(X)$ of a manifold with fibred corners X as the restriction of the groupoid of a decoupage in which X is embedded. One recovers the groupoid of $b$-calculus by Monthubert in the special case of a trivial fibration. The section ends with the proof of amenability and longitudinal smoothness of the groupoid $\Gamma_\phi(X)$.

Pseudodifferential calculus on manifolds with fibred corners is studied in the fourth section. To any manifold with fibred corners X we associate a pseudo\-differential calculus with compact support $\Psi_c^{\infty}(\Gamma_\phi(X))$ and we show it coincides with Melrose's $\phi$-calculus. We then introduce an extended calculus $\Psi^{\infty}(\Gamma_\phi(X)) = \Psi_c^{\infty}(\Gamma_\phi(X)) + \mathcal{S}_\psi(\Gamma_\phi(X))$ derived from a polynomial length function $\psi$ and we show this last algebra is stable under holomorphic calculus and includes the operators of the extended $\phi$-calculus of Mazzeo and Melrose.

Finally in the fifth section the groupoid of a fibred decoupage is shown to be the holonomy groupoid of the singular foliation defined by the manifold with fibred corners. The groupoid described in this article is therefore an explicit example of a singular leaf space in the sense of Androulidakis and Skandalis.


\tableofcontents

\section{Preliminaries on Lie groupoids}

The elements of index theory recalled hereafter mainly come from work by Connes \cite{Coinc}, Monthubert-Pierrot \cite{MP} et Nistor-Weinstein-Xu \cite{NWX}. The nice exposition in \cite{Ca3} is used as a guideline.

\subsection{Lie groupoids - Algebroids}

\subsubsection{Groupoids}
\begin{definition}
A $\textsl{groupoid}$ is a small category in which all morphisms are invertible. It is composed of a set of objects (or units) $\go$ and a set of arrows $\gr$ with two source and target applications $s,r:\gr \rightarrow \go$ and a composition law $m:\gr^{(2)}\rightarrow \gr$ associative on the set of composable arrows $\gr^{(2)}=\{ (\gamma,\eta)\in \gr \times \gr : s(\gamma)=r(\eta)\}$.
\end{definition}

It is also assumed there exists a unit map $u:\go \rightarrow \gr$ such that $r\circ u=s\circ u = \id$ and an inverse involutive map $i:\gr \rightarrow \gr$ such that $s\circ i=r$.

Finally by denoting $m(\gamma,\eta)=\gamma\cdot\eta$, every element $\gamma$ of $\gr$ satisfy the relations $\gamma\cdot\gamma^{-1}=u(r(\gamma))$, $\gamma^{-1}\cdot\gamma=u(s(\gamma))$ and $r(\gamma)\cdot\gamma=\gamma\cdot s(\gamma)=\gamma$.

A Lie groupoid is a groupoid for which $\gr$ and $\go$ are smooth manifolds, all maps above are smooth and $s,r$ are submersions.

It is classicaly denoted $\gr_A=s^{-1}(A)$, $\gr^B=r^{-1}(B)$ and $\gr_A^B=\gr_A\cap\gr^B$. 

\smallskip
\subsubsection{Algebroids}
\begin{definition}
Let M be a smooth manifold. A Lie algebroid on M is given by a vector bundle $A\rightarrow M$, a bracket $[ .\;,\;.]:\Gamma(A)\times\Gamma(A)\rightarrow \Gamma(A)$ over the module $\Gamma(A)$ of sections of $A$ and by a bundle morphism $p:A\rightarrow TM$ called \textsl{anchor map}, such that:
\begin{enumerate}
\item $[ .\;,\;.]$ is $\Rr$-bilinear, antisymetric et satisfies Jacobi identity,
\item  $[V,fW] = f[V,W] + p(V)(f)W$ for all $V,W\in\Gamma(A)$ and $f\in \ci(M)$
\item  $p([V,W]) = [p(V),p(W)]$ for all $V,W\in\Gamma(A)$.
\end{enumerate}
\end{definition}

Every Lie groupoid $\gr$ defines a Lie algebroid $\A\gr$ by the normal bundle of the inclusion $\go\subset \gr$. Sections of $\A\gr$ are in a bijective correspondance with vector fields on $\gr$ which are $s$-vertical and right-invariant, it therefore induces a Lie algebra structure on $\Gamma(\A\gr)$.

\begin{remark}
\label{algfol}
Every Lie algebroid defines a foliation by the image of its anchor map $p(\ci_c(M,A))$.
In particular every Lie groupoid defines a foliation.
\end{remark}

\subsection{$\gr$-pseudodifferential calculus}
A $\gr$-pseudodifferential operator is a differentiable family of pseudodifferential operators $\{P_x \}_{x \in \go}$ acting on $\ci_c(\gr_x)$ such that for $\gamma \in \gr$ and
$U_{\gamma}:\ci_c(\gr_{s(\gamma)}) \rightarrow \ci_c(\gr_{r(\gamma)}) $
the induced operator, the following $\gr$-invariance condition is verified:
$$ P_{r(\gamma)} \circ U_{\gamma}= U_{\gamma} \circ P_{s(\gamma)}.$$
Operators can act more generally on sections of a vector bundle $E \rightarrow \go $.
The exact condition of differentiability is defined in \cite{NWX}.

\noindent Only uniformly supported operators will be considered here, let us briefly recall this notion.
Let $P = (P_x, x \in \go)$ be a $\gr$-operator and denote $k_x$ the Schwartz kernel of $P_x$. 
The support and reduced support of P are defined by
$$supp\, P:=\overline{\cup_xsupp \,k_x}$$
$$supp_{\mu}P:=\mu_1(supp\,P)$$ where $\mu_1(g',g)=g'g^{-1}$.
We say that $P$ is uniformly supported if $supp_{\mu}P$ is compact in $\gr$.

In the following we denote $\Psi_c^m(\gr,E)$ the space of $\gr$-pseudodifferential operators uniformly supported. We denote also
$$\Psi_c^{\infty}(\gr,E)=\bigcup_m \Psi_c^m(\gr,E) \text{ and } \Psi_c^{-\infty}(\gr,E)=\bigcap_m \Psi_c^m(\gr,E).$$

\begin{remark}
The reduced support condition is justified by the fact that $\Psi_c^{-\infty}(\gr,E)$ is identified with
$\ci_c(\gr, End(E))$ using the Schwartz kernel theorem (\cite {NWX}).
\end{remark}
\begin{remark}
$\Psi_c^{\infty}(\gr,E)$ is a filtered algebra, $\textit{ie}$ $\Psi_c^{m}(\gr,E)\Psi_c^{m'}(\gr,E)\subset \Psi_c^{m+m'}(\gr,E).$ In particular, $\Psi_c^{-\infty}(\gr,E)$ is a two-sided ideal.
\end{remark}

The definition above is equivalent to that of \cite{Mont}, definition 1.1, where the space of pseudodifferential kernels on a longitudinally smooth groupoid $\gr$ is defined as the space $I^\infty(\gr,\go,\Omega^{\frac{1}{2}})$ of distributional sections $K$ on $\gr$ with values in half-densities $\Omega^{\frac{1}{2}}$ which are smooth outside $\go$ and given by oscillatory integral in a neighborhood of $\go$ :
$$ K(\gamma)=(2\pi)^{-n}\int_{A^*\gr_{r(\gamma)}}\exp{\left ( i<\phi(\gamma),\xi> \right )}a(\gamma,\xi)d\xi, $$
where $a$ is a polyhomogeneous symbol of any order with values in $\Omega^{\frac{1}{2}}$.

\smallskip

The notion of principal symbol extends easily to $\Psi_c^{\infty}(\gr,E)$. Denote by 
$\pi: A^*\gr \rightarrow \go$ the projection. If $P=(P_x,x\in \go)\in \Psi_c^{m}(\gr,E,F)$ is a pseudodifferential operator of order $m$ on $\gr$, 
the principal symbol $\sigma_m(P_x)$ of $P_x$ is a $\ci$ section of the vector bundle 
$End(\pi_x^*r^*E, \pi_x^*r^*F)$ over $T^*\gr_x$, such that the morphism defined over each fiber is homogeneous of degree $m$.
The existence of invariant connections allows the definition of exponentiation on the Lie algebroid $A^*\gr$ and provides a section $\sigma_m(P)$ of $End(\pi^*E, \pi^*F)$ over $A^*\gr$ which satisfies
\begin{equation}
\sigma_m(P)(\xi)=\sigma_m(P_x)(\xi)\in End(E_x,F_x) \text{ if } \xi \in A^*_x\gr
\end{equation}
modulo the space of symbols of order $m-1$. Terms of order $m$ of $\sigma_m(P)$ are invariant under a different choice of connection and the equation above induces a unique surjective linear map
\begin{equation}
\sigma_m:\Psi_c^{m}(\gr,E)\rightarrow \sw^m(A^*\gr,End(E,F)),
\end{equation}
with kernel $\Psi_c^{m-1}(\gr,E)$ (\cite{NWX}, proposition 2) where $\sw^m(A^*\gr,End(E,F))$ denotes the sections of the fiber $End(\pi^*E,\pi^*F)$ over $A^*\gr$ homogeneous of degree $m$ at each fiber.

\subsection{Index theory on singular spaces}
In the classical case of a compact manifold M without boundary, recall that an elliptic pseudodifferential operator D has a kernel and cokernel of finite dimension. The Fredholm index of D is the integer:
$$\ind D = \dim \ker D - \dim \coker D $$

This index is stable for any compact perturbation of D and its value only depends on topological data.
In the 60's  Atiyah and Singer introduced fundamental constructions in $K$-theory which allowed the understanding of the application ${D \rightarrow \ind D}$ through the group morphism
$$\ind_a: K^0(T^*M)\rightarrow\mathbb{Z},$$
called the \textsl{analytical index} of M.

More precisely, if $Ell(M)$ denotes the set of pseudodifferential operators on M, we have the following commutative diagram:
$$
\xymatrix{
Ell(M) \ar[r]^{\ind}\ar[d]_{\symb} & \mathbb{Z} \\
K^0(T^*M) \ar[ur]_{\ind_a} &
}
$$
where $Ell(M) \stackrel{\symb}{\rightarrow} K^0(T^*M)$ is the surjective application which maps an operator to the class of its principal symbol in $K^0(T^*M)$.
Atiyah and Singer then defined a topological index $\ind_t : K^0(T^*M) \rightarrow \mathbb{Z}$ with characteristic $K$-theoretical properties, showed that a unique morphism can check those properties and that the analytical index do satisfy them. The identity of the analytic and topological index is the celebrated Atiyah-Singer theorem \cite{AS0,AS,AS3}.

In the case of singular spaces such as the space of leaves of a regular foliation, orbifolds or manifolds with corners, similar constructions can be obtained by introducing \textsl{groupoids} \cite{Coinc,CS}.
A groupoid is a small category in which all morphisms are invertible. Groupoids
generalize the concepts of spaces, groups and equivalence relations. A groupoid is said to be \textsl{Lie} if the sets involved are smooth manifolds and morphisms are differential maps.
A pseudodifferential calculus was developed for Lie groupoids \cite{Coinc, MP, NWX} and more generally for longitudinally smooth groupoids such as \textsl{continuous family groupoids} \cite{LMN}.
 
The analytic index of elliptic operators on these groupoids is a morphism:
$$\ind_a :  K^0(A^*\gr) \rightarrow K_0(\cgr)$$
where $A^*\gr$ is the algebroid of the groupoid $\gr$ and $\cgr$ is the reduced $C^*$-algebra of $\gr$ which plays the role of the algebra of continuous functions on the singular space represented by $\gr$ \cite{Ren}.
Indices are thus generally not integers but elements of the $K$-theory group  $K_0(\cgr)$.
However in the case of a compact manifold M without boundary  they are classical indices as $\gr = M \times M$, $A^*\gr = T^*M$ and $K_0(\cgr) = K_0(\mathscr{K}) = \mathbb{Z}$. 

A fundamental property of the analytic index is that it can be factorized through the principal symbol map, 
there is a commutative diagram :
\[
\xymatrix{
Ell(\gr) \ar[r]^{ind}\ar[d]_{\sigma}& K_0(\cg)\ar[d]^j \\
K^0(A^*\gr) \ar[r]_{ind_a} & K_0(C^*_r(\gr)).
}.
\]

\noindent Indeed $ind_a$ is the boundary morphism associated with the short exact sequence of $C^{\ast}$-algebras (\cite{Coinc,CS,MP,NWX})
\begin{equation}
0\rightarrow \cgr\longrightarrow \overline{\Psi^0(\gr)}\stackrel{\sigma}{\longrightarrow}C_0(S^*\gr)\rightarrow 0
\end{equation}
where $\overline{\Psi^0(\gr)}$ is a $C^*-$completion of $\Psi_c^0(\gr)$, $S^*\gr$ the sphere bundle of  $A^*\gr$ and $\sigma$ the extension of the principal symbol.

\subsection{Normal cone deformation}
We recall in this subsection some basic properties of normal cone deformation groupoids, introduced by Hilsum and Skandalis \cite{HS}. These important objects in noncommutative geometry especially allow the construction of elements in $KK$-theory and the formulation of index problems (tangent groupoid of Connes \cite{Concg}). They are central in the description of the groupoid of manifolds with fibred corners.

\smallskip
\subsubsection{Definition}
Let $M$ be a smooth manifold and $N$ be a smooth submanifold of $M$. Denote $\mathscr{N}_N^M\rightarrow M$ the normal vector bundle to the inclusion $N\subset M$ : $\mathscr{N}_N^M\rightarrow M = T_NM/TN$.
The normal cone deformation of the inclusion $N\subset M$ is defined as a set by :
$$ D_N^M = \mathscr{N}_N^M \times 0 \bigsqcup M\times \Rr^*.$$

More generally, if $G_1$ and $G_2$ are two Lie groupoids with respective Lie algebroids $AG_1$ and $AG_2$ satisfying $AG_1\subset AG_2$, the normal cone deformation groupoid (\cite{HS}) can be defined for the induced immersion  $\varphi:G_1\rightarrow G_2$ by :
$$ D_\varphi = G_1\times_{s_1} \mathscr{N} \times 0 \bigsqcup G_2\times \Rr^*.$$

The groupoid $G_1$ acts on the normal $\mathscr{N}=AG_2/AG_1$ by its holonomy component $h_{xy}:\mathscr{N}_x\rightarrow \mathscr{N}_y$ :
$$(\gamma,v)=((x,y,h),(y,\xi))\mapsto \gamma\cdot v = (x,h^{-1}(\xi))$$

Thus $G_1\times_{s_V} \mathscr{N} = \{(\gamma,v) \in G_1\times \mathscr{N}, v\in \mathscr{N}_{s_1(\gamma)}\}$ is given a composition law by :
\begin{equation}
(\gamma_1,v)\cdot (\gamma_2,w)=(\gamma_1\gamma_2,v+\gamma_1\cdot w)
\end{equation}
when $t(\gamma_2)=s(\gamma_1)$, with inverse $(\gamma,v)^{-1}=(\gamma^{-1},-\gamma^{-1}\cdot v)$.

\noindent $G_2\times \Rr^*$ is naturally a Lie groupoid with point composition $(\gamma_1,t)\cdot (\gamma_2,t)=(\gamma_1\cdot\gamma_2,t)$.

\noindent $D_\varphi$ is then the union of the Lie groupoids $G_1\times_{s_V} \mathscr{N}$ and $G_2\times \Rr^*$.

\smallskip

\subsubsection{Differential structure}
\label{cone_diff}
In the simple case of the canonical inclusion $\Rr^p\times\{0\}\subset \Rp\times\Rq=\Rn$, where $q=n-p$, the normal cone deformation $D_p^n$ is the set $\Rp\times\Rq\times\Rr$ with the $\ci$ structure induced by the bijection $\Theta:\Rp\times\Rq\times\Rr\rightarrow D_p^n$:
\[
\Theta(x,\xi,t) = \left \{ \begin{array}{cc}
          (x,\xi,0)  & \text{ if } t=0 \\      
          (x,t\xi,t) & \text{ if } t\neq 0  
\end{array} \right. \]

In the local case of an open set $U\subset\Rn$ and a submanifold $V=U\cap(\Rp\times\{0\})$, $D_V^U = V\times \Rq\times\{0\}\bigsqcup U\times\Rr^*$ is an open subset of $D_p^n$ provided with the above structure as
$\Theta^{-1}(D_V^U)$ is the set
\begin{equation}
\label{omegauv}
\Omega_V^U = \{ (x,\xi,t) \in \Rp\times\Rq\times\Rr \; , \; (x,t\xi)\in U \}
\end{equation}
which is an open set of $\Rp\times\Rq\times\Rr$ and therefore a smooth manifold.

In the general case of manifolds suppose $M$ of dimension $n$ and $N$ of dimension $p$. The differential structure on $D_N^M$ can be described locally from the open sets $D_p^n$ with adapted charts.

\begin{definition}[\cite{Ca3}]
A local chart $(\mathcal{U},\phi)$ of $M$ is said to be a $N$-slice if
\begin{enumerate}
\item $\phi:\mathcal{U}\stackrel{\simeq}{\rightarrow} U\subset \Rp\times\Rq$
\item if $\mathcal{U}\cap N = \mathcal{V}$, $\mathcal{V} = \phi^{-1}(V)$ (where $V=U\cap(\Rp\times\{0\})$ as above)
\end{enumerate}
\end{definition}

Let $(\mathcal{U},\phi)$ be a $N$-slice. When $x\in\mathcal{V}$ we have $\phi(x)\in \Rp\times\{0\}$ : we denote by $\phi_1$ the component of $\phi$ on $\Rp$ such that $\phi(x)=(\phi_1(x),0)$ and $d_n\phi_v : \mathcal{N}_v\rightarrow \Rq$ the normal component of the differential $d\phi_v$ for $x\in\mathcal{V}$. 
Let $\tilde{\phi}:D_\mathcal{U}^\mathcal{V} \rightarrow D_U^V$ be the map defined by :
\[
\tilde{\phi}:\left \{  \begin{array}{cccc}
             \tilde{\phi}(v,\xi,0) & = & (\phi_1(v),d_n\phi_v(\xi),0) &\\      
             \tilde{\phi}(u,t)  & =  & (\phi(u),t) & \text{ if } t\neq 0
\end{array} \right. \]

Then the map $\varphi = \Theta^{-1}\circ\tilde{\phi}$ obtained as the composition
$$ D_\mathcal{V}^\mathcal{U} \stackrel{\tilde{\phi}}{\longrightarrow} D_V^U \stackrel{\Theta^{-1}}{\longrightarrow}\Omega_V^U $$
is a diffeomorphism of $D_\mathcal{V}^\mathcal{U}$ on the open set $\Omega_V^U$ of $\Rp\times\Rq\times\Rr$ defined in the same way as $(\ref{omegauv})$.

The global differential structure of $D_N^M$ is then described by the following proposition :
\begin{proposition}[\cite{Ca3}]
Let $\{(\mathcal{U}_{\alpha},\phi_{\alpha})\}$ be a $\ci$ atlas of M composed of N-slices. 
Then $\{(D_{\mathcal{V}_{\alpha}}^{\mathcal{U}_{\alpha}},\varphi_{\alpha})\}$ is a $\ci$ atlas of $D_N^M$.
\end{proposition}

A change of atlas is described by the following result (see \cite{Ca3}) :
Let $U\subset \Rp\times\Rq$ and $U^\prime\subset \Rp\times\Rq$ be open sets and $F:U\rightarrow U^\prime$ a $\ci$ diffeomorphism decomposed as a $\Rp\times\{0\}$-slice $F=(F_1,F_2)$ such that $F_2(x,0)=0$.
Then the map $\tilde{F}:\Omega_V^U\rightarrow \Omega_{V^\prime}^{U^\prime}$ defined by :
\[
\tilde{F}(x,\xi,t) = \left \{ \begin{array}{cc}
          \frac{\partial F}{\partial\xi}(x,0)\cdot\xi  & \text{ if } t=0 \\      
          \frac{1}{t}F(x,t\xi) & \text{ if } t\neq 0  
\end{array} \right. \]

is a $\ci$ diffeomorphism from $\Omega_V^U$ to $ \Omega_{V^\prime}^{U^\prime}$.

\section{Manifolds with fibred corners}

\subsection{Manifolds with corners}
A manifold with corners $X$ is a topological space where any point $p$ has a neighbourhood diffeomorphic to $\Rnk$, with 0 as image of $p$. $k$ is called the \textsl{codimension} of $p$. The connected components of the set of points with codimension $k$ are called \textsl{open faces} of codimension $k$. Their closure are the \textsl{faces} of $X$. A face of codimension 1 is called an \textsl{hyperface} of $X$. The boundary of X is denoted by $\partial X$, it is the union of faces with codimension $k>0$.

\subsection{Manifolds with embedded corners}
Manifolds with embedded corners have been studied by Melrose to understand and extend Atiyah-Patodi-Singer index theorem (see \cite{Mel}). An equivalent description in terms of decoupages has been proposed by Monthubert in \cite{Mont}. The two definitions are recalled below (section \ref{equar} for the definition with decoupages).

\begin{definition}[Manifold with embedded corners, \cite{Mel}]
A manifold with embedded corners $X$ is a manifold with corners endowed with a subalgebra $\ci(X)$ satisfying the following conditions:
\begin{itemize}
	\item[$\cdot$] there exists a manifold $\mX$ and a map $j:X\rightarrow \mX$ such that $$\ci(X)=j^*\ci(\mX),$$ 
	\item[$\cdot$] there exists a finite family of functions $\rho_i \in \ci(\mX)$ such that $$j(X)=\{y\in\mX,\forall i \in I,\rho_i(y)\geq 0\},$$
	\item[$\cdot$] for any $J\subset I$ and any point of $\mX$ where the $(\rho_j)_{j\in J}$ simultaneously vanish, the differentials d$\rho_j$ are independent.
\end{itemize}
The functions $\rho_i$ are called \textsl{defining functions} of the hyperfaces.
\end{definition}

\begin{remark}
Manifolds with embedded corners are manifolds with corners : the number of functions $\rho_j$ which vanish at a point $p$ determine the codimension of $p$ since the $\rho_j$, which are independent, can be used as local coordinates and induce a diffeomorphism from a neighbourhood of $p$ to some $\Rnk$.
\end{remark}
\begin{remark}
A definition function $\rho$ of an hyperface F induces a trivialization of the normal bundle on its interior $F^\circ$.
\end{remark}

\begin{examples}
The edge, the square, cubes of dimension $n$ are manifolds with embedded corners.
The simplex $\Sigma_n = \{ (t_0,\dots,t_n)\in \Rr^{n+1}_+, t_0+\dots+t_n = 1 \}$ is a manifold with embedded corners.
The drop with one corner is not a manifold with embedded corners, contrary to the drop with two corners.
\end{examples}

\begin{figure}[hbtp]
\begin{minipage}[b]{0.3\linewidth}
\centering
	\begin{pspicture}(-1,-1)(2,2)
    \psarc(1,1){1}{180}{90}
    \psline(1,2)(0,2)(0,1)
	\end{pspicture}
\end{minipage}
\hspace{1cm}
\begin{minipage}[b]{0.3\linewidth}
\centering
	\begin{pspicture}(0,-2)(3,2)
	\psbezier[linewidth=1pt](0,0)(1.5,2)(3,0)
    \psbezier[linewidth=1pt](0,0)(1.5,-2)(3,0)
	\end{pspicture}
\end{minipage}\vspace*{\fill}
\label{coins}
\caption{Drops with one and two corners.}
\end{figure}
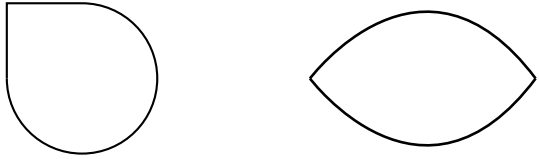

\subsection{Manifolds with iterated fibred corners}
The study of pseudo-stratified manifolds motivates the introduction of a fibred structure on manifolds with embedded corners (\cite{DLR}). Indeed the existence of a desingularization process shows that any pseudo-stratified manifold is a quotient space of a manifold with fibred corners (\cite{DL},\cite{DLR},\cite{ALMP}). 
For such a desingularization each hyperface $F_i$ is the total space of a fibration $\phi_i:F_i\rightarrow Y_i$ where $Y_i$ is also a manifold with corners such that the family of fibrations $\phi=(\phi_1,\dots,\phi_k)$ satisfy the properties of an \textsl{iterated fibred structure} (\cite{DLR}), ie :
\begin{itemize}
	\item[$\cdot$] $\forall I\subset \{1,\dots,k\} $ such that $\bigcap_{i\in I} F_i\neq\emptyset$, the set  $\{F_i,i\in I\}$ is totally ordered,
	\item[$\cdot$] if $F_i<F_j$, then $F_i\cap F_j\neq\emptyset$, $\phi_i:F_i\cap F_j\rightarrow Y_i$ is a surjective submersion and $Y_{ji}\doteq\phi_j(F_i\cap F_j)\subset Y_j$ is an hyperface of the manifold with corners $Y_j$. Moreover there exists a surjective submersion $\phi_{ji}:Y_{ji}\rightarrow Y_i$ which satisfies  $\phi_{ji}\circ\phi_j=\phi_i$ on $F_i\cap F_j$. 
	\item[$\cdot$] the hyperfaces of $Y_j$ are exactly the $Y_{ji}$ with $F_i<F_j$. In particular if $F_i$ is minimal $Y_i$ is a manifold without boundary.
\end{itemize}  

We call \textsl{manifold with iterated fibred corners} a manifold with embedded corners endowed with an iterated fibred structure.

\subsection{Decoupages}

We recall the definition of a \textsl{decoupage} which is a useful notion to associate to any manifold with embedded corner a ``puff" Lie groupoid \cite{Mont}. Generalizing this approach we then introduce an original definition of \textsl{fibred decoupage} suited to the case of fibrations and allowing the direct construction of a puff groupoid for manifolds with fibred corners \cite{Gui}. Before let us recall a few facts about transversality.

\begin{definition}
Let $X_1,\dots,X_n$ be a family of smooth manifolds and for each $i$ let $\phi_i:X_i\rightarrow Y$ be a smooth map. If $(x_i)_{1\leq i\leq n}\in \prod_{1\leq i\leq n}X_i$ is a family such that $\phi_1(x_1)=\dots = \phi_n(x_n)$, the family $(\phi_i)_{1\leq i\leq n}, \phi_i:X_i\rightarrow Y$ is said to be \textsl{transverse} if the orthogonals (in $T^*Y$) of the spaces d$\phi_i(TX_i)$ are in direct sum. 
\end{definition}

Under this condition the fibered product of the $X_i$ over Y is a smooth submanifold of $\prod_{1\leq i\leq n}X_i$. 

A family of submanifolds will be called transverse if the inclusions of these manifolds are transverse; a smooth map is said to be transverse if it is transverse to the inclusion of this submanifold.

\smallskip

\begin{definition}[Decoupage, \cite{Mont} 2.4]
\label{equar}
A \textsl{decoupage} is given by a manifold $M$ and a finite family $(V_i)_{i\in I}$ of submanifolds of codimension 1 such that $\forall J\subset I$ the family of inclusions of the $(V_j)_{j\in J}$ is transverse.
\end{definition}

The decoupage is said to be \textsl{oriented} if every submanifold $V_i$ is transversally oriented. Moreover if each $V_i$ splits $M$ in two parts $M^+_i$ and $M^-_i$ the intersection of the $M^+_i$ is called \textsl{positive part}. 
 
The equivalent definition of a manifold with embedded corners as a decoupage is then:
\begin{definition}[Manifold with embedded corners]
A manifold with embedded corners $X$ is the positive part of an oriented decoupage $(M,(V_i)_{i\in I})$.
$X$ is the positive part of this decoupage and $(M,(V_i)_{i\in I})$ is said to be an \textsl{extension} of $X$.
\end{definition}

\begin{definition}[Fibred decoupage, \cite{Gui} 2.3]
A \textsl{fibred decoupage} is given by a decoupage $(M,(V_i)_{i\in I})$ and a family of fibrations $\phi_i:V_i\rightarrow Y_i$ such that $\forall J\subset I$ the family of inclusions of the $(V_j\times_{Y_j} V_j)_{j\in J}$ in $M^2$ is transverse.
\end{definition}

\subsection{Manifolds with fibred corners}
\begin{definition}[\cite{Gui} 2.4]
Let $X$ be a manifold with embedded corners and $(M,(V_i)_{i\in I})$ an extension of $X$.
We call $X$ a \textsl{manifolds with fibred corners} if there exists a family of fibrations $\phi_i:V_i\rightarrow Y_i$ such that $(M,(V_i,\phi_i)_{i\in I})$ is a fibred decoupage.
\end{definition}

\begin{example}
A family of disjoint fibred submanifolds of codimension 1 is always transverse. In particular manifolds with fibred boundaries, the objects of Melrose's $\phi$-calculus, are manifolds with fibred corners.
\end{example}

\begin{example}
\label{ex1}
$X=\Rr_+^3$. Let $M=\Rr^3$, $\pi_x$, $\pi_y$ and $\pi_z$ be the projections on the canonical basis, $V_x=\ker\pi_x$, $V_y=\ker\pi_y$, $V_z=\ker\pi_z$. 
Let $\phi_x:V_x\rightarrow \{0\}$, $\phi_y:V_y\rightarrow \ker\pi_x\cap  \ker\pi_y$, $\phi_z:V_z\rightarrow \ker\pi_y\cap  \ker\pi_z$.
Then $(M,(V_x,\phi_x),(V_y,\phi_y),(V_z,\phi_z))$ is a fibred decoupage.
\end{example}

\begin{example}
Same notations as previous example but with $\phi'_x:V_x\rightarrow \{0\}$, $\phi'_y:V_y\rightarrow \ker\pi_y\cap  \ker\pi_z$, $\phi'_z:V_z\rightarrow \ker\pi_x\cap  \ker\pi_z$.
Then $(M,(V_x,\phi'_x),(V_y,\phi'_y),(V_z,\phi'_z))$ is not a fibred decoupage.
\end{example}

\begin{example}
\label{ex3}
The example \ref{ex1} is a fibred decoupage but is not a iterated fibred structure.
Indeed $\phi_z(V_y\cap V_z)=V_z\cap\ker\pi_y$, $Y_{yz}=\phi_y(V_y\cap V_z)= \{0\}$ and there is no surjection of  $Y_{yz}$ on $Y_z=\ker\pi_y\cap\ker\pi_z$ as $\dim Y_{yz}=0<1=\dim Y_z$.
\end{example}

\begin{remark}
The definition used for manifolds with fibred corners ensures longitudinal smoothness of the puff groupoid (see Proposition \ref{long_eclat}). The results of this chapter can therefore be applied to fibrations coming from other structures than  iterated fibred structures, since any fibred decoupage does not necessarily come from an iterated fibred structures (eg \ref{ex3}).

However, it is perfectly possible to limit the study to the category of \textsl{iterated fibred decoupages}, ie the decoupages which fibrations satisfy the properties of an iterated structure.
\end{remark}

\section{The groupoid of manifolds with fibred corners}

\subsection{The groupoid of a codimension 1 foliated submanifold}
In this subsection we construct a ``puff groupoid" $\grVF$ for any foliated codimension 1 submanifold $(V,\mathcal{F})$ of a manifold M. This groupoid is described as a gluing between the holonomy groupoid of $M\setminus V$ and a deformation groupoid $\D_\varphi = D_\varphi\rtimes\Rr_+^*$ built from the normal cone deformation $D_\varphi$ of the groupoid immersion $\mathscr{H}ol(\mathcal{F})\rightarrow V\times V$:
$$\gr V_\mathcal{F} = \D_\varphi \bigcup_\Psi \mathcal{H}ol(M\setminus V).$$

Gluing is preformed through a map $\Psi$ and one has the following commutative diagram
\[
\xymatrix{
V\times V\times\Rr_+^*\times \Rr^* \ar[r]^{\Psi}\ar[d]_{i_A}& \mathcal{H}ol(M\setminus V) \ar[d]^{i_{M\setminus V}} \\
\D_\varphi \ar[r]_{i_\D} & \gr V_\mathcal{F}
}
\]

$\grVF$ encodes the space of leaves of the singular foliation $\F$ defined by vector fields on $M$ tangent to $\mathcal{F}$ on $V$. This property will be precised in the section devoted to singular foliations : we will show  $\grVF$ is nothing but the holonomy groupoid $\HolM$ of $\F$ (proposition \ref{gr_holonomy}). 

\smallskip

\subsubsection{Notations}
\label{notations}
Let $M$ be a smooth manifold and $V\subset M$ a connected submanifold of codimension 1 transversally oriented. $V$ is supposed to be foliated by a regular foliation $\mathcal{F}$ defined by an integrable vector subspace  $TF\subset TV$, i.e. $\ci(V,TF)$ is a Lie subalgebra of $\ci(V,TV)$. Let denote $\mathscr{N} = TV/TF$ and suppose $M$ is given a connexion $w$ which restriction to $V$ provides a decomposition of the tangent bundle  $TV = TF \oplus \mathscr{N}$.

Let $(N,\pi,V)$ be a tubular neighbourhood of V in M : N is an open set of M including V, $\pi:N\rightarrow V$ a vector bundle with $\Rr$ type fibre. 
Let $\{f_i:N_i=\pi^{-1}(V_i)\rightarrow V_i\times\Rr\}_{i\in I}$ be a local trivialization of N, which components are denoted $f_i = (\pi_i,\rho_i)_{i\in I}$. 

\smallskip

\subsubsection{The deformation groupoid $\D_\varphi = D_\varphi\rtimes\Rr_+^* $}
\label{d_varphi}

\paragraph{Definition} 
Let $G_1$ be the holonomy groupoid \cite{Cosf,Concg} associated with $\mathcal{F}$ and ${\varphi:G_1 \rightarrow V\times V }$ the immersion induced by the inclusion  $TF \subset TV$. The normal cone deformation groupoid of $\varphi$ is denoted $D_\varphi$. 

$D_\varphi$ induces on its space of units $V\times\Rr$ a singular foliation $\mathcal{F}_{\text{adp}}$ defined by $\mathcal{F}\times \{0\}$ and $TV\times \{t\}$ for $t\neq 0$ resulting from the integration of the partial adiabatic Lie algebroid $\mathscr{A}_{\text{adp}}$ described in (\cite{Deb}, Example 4.5): 

\begin{center}
\begin{math}
   \begin{array}{cccc}
     p_{\text{adp}}: & \mathscr{A}_{\text{adp}} =TF\oplus\mathscr{N}\times\Rr &  
                                                           \rightarrow & TV\times T\Rr  \\
      & (v_1,v_2,t) &  \mapsto & (v_1+tv_2,(t,0))
    \end{array}  
\end{math}
\end{center}

\noindent In particular the foliation of $\Rr$ transverse to $V$ induced by $D_\varphi$ is trivial in points as each $t\neq 0$ defines a leaf $V\times\{t\}$ of $V\times\Rr$.

But the foliation $\mathcal{F}_M$ defined by the transversally oriented faces $V_i$ of a manifold with fibred corners has a singular transverse structure of the form $V_i\times\{0\}$,$V_i\times\Rr^*_+$ and $V_i\times\Rr^*_-$. For that purpose we introduce the deformation groupoid $\D_\varphi$ obtained by the transverse action of $\Rr_+^*$ on $D_\varphi$ explicitely given by :
$$(\lambda,h)=(\lambda,(\gamma,t)) \mapsto \lambda\cdot \gamma = (\gamma,\lambda\cdot t)$$

The orbits of $\D_\varphi=D_\varphi\rtimes\Rr_+^*$ then coincide with the foliation $\mathcal{F}_M$. 

An expression of the manifold $\D_\varphi$ as a set is :
$$\D_\varphi=(G_1\times_{s_1} \mathscr{N} \times \Rr_+^*)\times\{0\} \bigsqcup V\times V\times \Rr^*_+\times\Rr^*\rightrightarrows V\times\Rr.$$

\paragraph{Composition law}
Let $p:V\times \Rr^*_+\times \Rr\rightarrow V$ be the canonical projection, $\hat{s}=p\circ s_\varphi$ and $\hat{r}=p\circ r_\varphi$. For a given element $g=(\gamma,\lambda,t)$ of $\D_\varphi$, $\hat{s}(g)$ and $\hat{r}(g)$ only depend on $\gamma\in G_1\times_{s_1}\mathscr{N} \bigsqcup V\times V$. To keep notations simple we denote $\hat{s}(\gamma)$ and $\hat{r}(\gamma)$ the corresponding elements.

The composition law on $\D_\varphi$ is then the product of the law composition on $D_\varphi$ as a normal cone groupoid and the law on $H=\Rr\rtimes\Rr_+^*$ where $\Rr_+^*$ acts on $\Rr$ by multiplication. If $g_1 = (\gamma_1,\lambda_1,t_1)$ and $g_2 = (\gamma_2,\lambda_2,t_2)$ are two elements of $\D_\varphi$ such that $\hat{s}(\gamma_2)=\hat{r}(\gamma_1)$ and $s_H(\lambda_2,t_2) = r_H(\lambda_1,t_1)$, we have :
\begin{equation}
\label{loi_compo}
  g_1\cdot g_2  =  (\gamma_1\cdot\gamma_2,\lambda_1 \cdot \lambda_2, t_1)
\end{equation}

\smallskip

\subsubsection{The puff groupoid $\grVF$}
We are ready to define the gluing function $\Psi$. 
Let $A \subset \D_\varphi$ be the restriction of $\D_\varphi$ to the product component $V\times V\times \Rr_+^*\times\Rr^*$ and set $A_i^j=\A\cap\D_i^j=V_i\times V_j\times \Rr_+^*\times\Rr^*$ where
$\D_i^j=\{(\gamma,\lambda,t)\in \D_\varphi \;/\; \hat{s}(\gamma)\in V_i,\; \hat{r}(\gamma)\in V_j \}.$ 

Now consider the maps  $\Psi_{ij}:A_i^j\rightarrow M^\circ\times M^\circ$ defined by
 \[ \Psi_{ij}(\gamma,\lambda,t) = 
          ( f_i^{-1}(\hat{s}(\gamma),t), f_j^{-1}(\hat{r}(\gamma),\lambda \cdot t)) 
 \]
 
The $A_i^j$ being disjoint we have $A=\sqcup_{i,j}A_i^j$ and one defines $\Psi:A\rightarrow  M^\circ\times M^\circ$ by $\Psi(a)=\Psi_{ij}(a)$ for $a\in A_i^j$. 
It is immediate to check that $\Psi$ preserves the transverse orientation of V in M : $\Psi(A_i^{j+})\subset (M^+\setminus V)^2$, $\Psi(A_i^{j-})\subset (M^-\setminus V)^2$. V being connected, the image of $\Psi$ is thus in the union of the connected components of $(M\setminus V)^2$, that is to say in $\mathcal{H}ol(M\setminus V)$.

$\gr V_\mathcal{F}$ is then defined as a topological space by the gluing of $\D_\varphi$ and $\mathcal{H}ol(M\setminus V)$ by the map $\Psi$:
$$ \gr V_\mathcal{F} = \D_\varphi \bigcup_\Psi \mathcal{H}ol(M\setminus V) $$

and one has the following commutative diagram :
\[
\xymatrix{
V\times V\times\Rr_+^*\times \Rr^* \ar[r]^{\Psi}\ar[d]_{i_A}& \mathcal{H}ol(M\setminus V) \ar[d]^{i_{M\setminus V}} \\
\D_\varphi \ar[r]_{i_\D} & \gr V_\mathcal{F}.
}
\]

\smallskip

\subsubsection{Differential structure of the puff groupoid $\grVF$}
\label{hol_diff}
Let $(\varphi_{ij})_{i,j\in I^2}$ be an atlas on $\D_\varphi$ which charts are slices of ${\D_i^j\cap (G_1\times \Rr_+^*)}$.
Set $p=\dim{G_1}$ and $q=\dim{\mathscr{N}}=2\cdot\dim{V}-p$.

Let $U_i^j$, W and Z be respective open sets of $\Rr^{p+q}, \Rr^p$ and $\Rr$ satisfying ${W=U_i^j\cap(\Rr^p\times\{0\})}$ and $\tilde{\varphi_{ij}}^{-1}(U_i^j\times\Rr^*)\subset V_i\times V_j\times\Rr^*$.

Topology on $\D_\varphi$ is generated (see section \ref{cone_diff}) by the open sets ${\Theta\circ\tilde{\varphi_{ij}}(\Omega_W^{U_i^j})\times\exp{Z}}$ where:
$$\Omega_W^{U_i^j} = \{ (x,\xi,t) \in \Rp\times\Rq\times\Rr \; , \; (x,t\xi)\in U_i^j \}.$$

In particular $\D_i^j$ is an open set of $\D_\varphi$ generated by the image of $\Omega_W^{V_i\times V_j}$ and $Z=\Rr$.
Similary $A_i^j$ is an open set of $\D_\varphi$ generated by the image of $Z=\Rr$ and $\Omega_*^{V_i\times V_j}=\Omega_W^{V_i\times V_j}\cap (\Rp\times\Rq\times\Rr^*)$.

One then checks that $s_\varphi\times r_\varphi$ is an open map on $A^2$ as the map $(t,\lambda)\mapsto (t,\lambda\cdot t)$ is a diffeomorphism of $\Rr^*\times\Rr_+^*$. 
The map $\Psi_{ij}:A_i^j\rightarrow M^\circ\times M^\circ$ is thus open as a composition of the open maps $f_i\times f_j$ and $s_\varphi\times r_\varphi$.

Therefore $i_\D:\D_\varphi\rightarrow\grVF$ defined by :
 \[ i_\D(g) = \left \{
                            \begin{array}{cc}
                             \Psi(g) & \text{ if } g\in A \\
                              g  & \text{ otherwise }
                            \end{array}
                    \right.
 \]
is an open map and an atlas $\mathcal{A} = \{(\Omega_i^j,\phi_{ij})_{i,j\in I^2}\}$ on $\grVF$ is obtained by setting $\Omega_i^j=i_\D(\D_i^j)$ and $\phi_{ij}=\varphi_{ij}\circ i_{\D}^{-1}$.

Let $\tilde{\mathcal{B}}$ be an atlas of the smooth manifold $M^{\circ}\times M^{\circ}$ and let $\mathcal{B} = \{(\Omega_\beta,\varphi_\beta)\}$ be the induced atlas on the manifold $\mathcal{H}ol(M\setminus V)$.

The differential structure on $\grVF$ is obtained as the product structure of $M^{\circ}\times M^{\circ}$ and the differential structure induced by the family $(\phi_{ij})$. More precisely we have the following proposition :
\begin{proposition}
$\mathcal{A}\cup\mathcal{B}$ is a smooth atlas on $\grVF$.
\end{proposition}
\begin{demo}
Let $(\Omega,\varphi)$ be a chart of $\mathcal{B}$.
We must prove that any chart of $\mathcal{A}$ meeting $\Omega$ is compatible with $\varphi$.
Let $i$ and $j$ be such that $\Omega_i^j\cap\Omega\neq\emptyset$ and define $\Phi_{ij}:\varphi(\Omega_i^j\cap\Omega)\rightarrow \phi_{ij}(\D_i^j)$ by $\Phi_{ij}=\phi_{ij}\circ\Psi_{ij}^{-1}\circ\varphi^{-1}$.
It is immediate to see that $\Omega_i^j\cap\Omega \subset \Psi_{ij}(A_i^j)=N^+_i\setminus V_i\times N^+_j\setminus V_j \bigsqcup N^-_i\setminus V_i\times N^-_j\setminus V_j$ and that 
 \[ 
\Psi_{ij}^{-1}: \begin{array}{ccc}
                   \Omega_i^j\cap\Omega & \rightarrow & A_i^j \\
                        (x,y) & \rightarrow &  \left ( \pi_i (x),\pi_j (y) , \rho_j(y) / \rho_i (x),\rho_i (x) \right )
                \end{array}
 \] 

\noindent is a diffeomorphism on its image.
It implies that $\Phi_{ij}$ is a $C^{\infty}$ clutching function and the compatibility of $\mathcal{A}$ and $\mathcal{B}$ follows.
\end{demo}

\begin{remark}
In the general case above $\D_\varphi$ only describes the differential structure of $\gr V_\mathcal{F}$ in a neighbourhood of V.
The description gets simpler if a trivialisation $(\pi,\rho)$ is given with $\pi$ and $\rho$ globally defined on  M. 
The map $\Psi$ is then bijective, the inverse map $\Psi^{-1}:\mathcal{H}ol(M\setminus V) \rightarrow V\times V\times \Rr_+^*$ being defined by : $$\Psi^{-1}(x,y)= \left( \pi(x),\pi(y),\rho(y)/\rho(x),\rho(x) \right).$$
Thus $i_A\circ\Psi^{-1}$ is injective and $i_{D}$ is a diffeomorphism of $D_\varphi$ on $\im i_{D} = \gr V_\mathcal{F}$.
\end{remark}

\smallskip

\subsubsection{Composition law}
The groupoid $G_1=\mathscr{H}ol(\F)$ acts on the normal bundle $\mathscr{N}$ by its holonomy component $h_{xy}:\mathscr{N}_x\rightarrow \mathscr{N}_y$ :
$$(\gamma,v)=((x,y,h),(y,\xi))\mapsto \gamma\cdot v = (x,h^{-1}(\xi))$$

Therefore $G_1\times_{s_V} \mathscr{N} \times \Rr_+^* = \{(\gamma,v,\lambda) \in G_1\times \mathscr{N} \times \Rr_+^*, v\in \mathscr{N}_{s_1(\gamma)}\}$ is endowed with a composition law :
\begin{equation}
\label{loicompo}
(\gamma_1,v,\lambda_1)\cdot (\gamma_2,w,\lambda_2)=(\gamma_1\gamma_2,v+\gamma_1\cdot w,\lambda_1\cdot \lambda_2)
\end{equation}
when $t(\gamma_2)=s(\gamma_1)$, with inverse $(\gamma,v,\lambda)^{-1}=(\gamma^{-1},-\gamma^{-1}\cdot v,\lambda^{-1})$.

The composition law on $\gr V_\mathcal{F}$ is defined as the law (\ref{loicompo}) on $\mathscr{N}_{\mathcal{F}} = G_1\times_{s_1} \mathscr{N} \times \Rr_+^*$ and the canonical product law on $M^{\circ}_+\times M^{\circ}_+ \bigsqcup M^{\circ}_-\times M^{\circ}_-$.

Source and target maps of $\gr V_\mathcal{F}$ are defined on $\mathscr{N}_{\mathcal{F}_V}$ and its complementary by:
\[
s: \left \{ \begin{array}{lll}
           s(\gamma,\xi,\lambda) & = & s_1(\gamma)  \\    
           s(x,y)                & = & x 
\end{array} \right.
\]
and
\[
r: \left \{ \begin{array}{lll}
           r(\gamma,\xi,\lambda) & = & r_1(\gamma)  \\    
           r(x,y)                & = & y 
\end{array} \right. 
\]

\bigskip 
 
The groupoid structure of $\gr V_\mathcal{F}$ can be easily reformulated on the deformation groupoid $\D_\varphi$. In fact $\D_\varphi$ endowed with its composition law \ref{loi_compo} is a Lie groupoid isomorphic to a neighbourhood of $\mathscr{N}_{\mathcal{F}}$ in $\gr V_\mathcal{F}$. More precisely, let $i_{\D}$ be the inclusion of $\D_\varphi$ in $\gr V_\mathcal{F}$ given by the gluing $\Psi$. Then :

\begin{proposition}
\label{psi_G}
The Lie groupoids $\D_\varphi$ and $i_{\D}(\D_\varphi)\subset \gr V_\mathcal{F}$ are isomorphic.
\end{proposition}

\begin{demo}
$i_{\D}$ is by construction a diffeomorphism from $\D_\varphi$ to $i_{\D}(\D_\varphi)$.
We observe that $i_{\D}$ can be written under the form :
 \[ i_{\D} : \left \{
    \begin{array}{cll}
      g=(\gamma,\lambda,t)\in \D_i^j \mapsto &  ( f_i^{-1}\circ s_\varphi(g), f_j^{-1}\circ r_\varphi(g)) & \text{ if } t\neq 0 \\
                                      &  ( \gamma,\lambda, 0 )  & \text{ if } t=0
    \end{array}    \right.
 \]
Let $g_1=(\gamma_1,\lambda_1,t_1)$ and $g_2=(\gamma_2,\lambda_2,t_2)$ be two elements of $\D_\varphi$.

The relation $i_{\D}(g_1)\cdot i_{\D}(g_2) = i_{\D}(g_1 \cdot g_2)$ is trivial if $t_1=0$ since then $i_{\D} = \text{Id}$.

For $t_1 \neq 0$,  $g_1\cdot g_2 = (\gamma_1\cdot\gamma_2,\lambda_1 \cdot \lambda_2, t_1)$.
One computes :
 \[
        \begin{array}{lcl}
i_{\D}(g_1)\cdot i_{\D}(g_2) & = & ( f_i^{-1}\circ s_\varphi(g_1), f_j^{-1}\circ r_\varphi(g_1)) \cdot  ( f_i^{-1}\circ s_\varphi(g_2), f_j^{-1}\circ r_\varphi(g_2))   \\    
          & = & ( f_i^{-1}\circ s_\varphi(g_1), f_j^{-1}\circ r_\varphi(g_2)) \text{ with } t_2 = \lambda_1\cdot t_1 \\
          & = & ( f_i^{-1}\circ s_\varphi(g_1), f_j^{-1}\circ (\hat{r}(\gamma_2),\lambda_2\cdot \lambda_1 \cdot t_1) )  \\  
          & = & i_{\D}(g_1 \cdot g_2)  
        \end{array} 
\]

Thus $i_{\D}(g_1)\cdot i_{\D}(g_2) = i_{\D}(g_1 \cdot g_2)$ for any composable elements of $\D_\varphi$ and $i_{\D}$ is indeed a Lie groupoid morphism.

\end{demo}

\subsubsection{The case of fibrations}
Assume the foliation $\mathcal{F}$ of V is a fibration $\phi:V\rightarrow Y$.
The immersion $\varphi:G_1\rightarrow G_V$ is then an embedding with $G_1=V\times_YV$ and $G_V=V\times V$.
The holonomy of a fibration is trivial and $G_1$ trivially acts on $\mathscr{N}$ by:
$$(\gamma,v)=((x,x^\prime),(x^\prime,\xi))\mapsto \gamma\cdot v = (x,\xi).$$

Moreover for a local trivialization also trivializing the fibration $\phi:V\rightarrow Y$, $\mathscr{N}$ can be identified as the pushout of $TY$. The puff groupoid is then locally described as a set by:
$$\gr V_{\phi|N_i}=(V_i\times_YV_i\times_YTY\times\Rr_+^*)\bigsqcup (N^+_i\setminus V_i)\times (N^+_i\setminus V_i) \bigsqcup (N^-_i\setminus V_i)\times (N^-_i\setminus V_i).$$
\label{gr_loc}

The two limit cases of a coarse fibration $Y=\{pt\}$ and a ponctual fibration $Y=V$ correspond to the geometric situations of $b$-calculus and 0-calculus. The normal cone deformations $D_{\varphi}$ are then respectively isomorphic to the product groupoid $V\times V$ and to the adiabatic groupoid of V. The puff of the coarse fibration is the groupoid described in \cite{Mont} and \cite{NWX}. The definition of $\gr V_\mathcal{F}$ as a holonomy groupoid $\HolM$ of the foliation $\F$ (proposition \ref{gr_holonomy}) gives those groupoids a geometric meaning independent of any particular construction.

\begin{remark}
A similar construction allows the description of a manifold which interior is endowed with a coarser fibration $\Phi:M\rightarrow B$: if $G_\Phi=V\times_\Phi V$ denotes the holonomy groupoid of the fibration $\Phi$ restricted to the boundary, $\gr V_\phi$ is obtained as the gluing of $\D_\varphi$ and $\mathcal{H}ol(M,\Phi)=M\times_B M$, where $\D_\varphi=D_\varphi\rtimes\Rr_+^*$ is built from the normal cone deformation $\varphi:G_\phi\rightarrow G_\Phi$.
\end{remark}

\subsection{The groupoid of manifolds with fibred corners}

The groupoid of a fibred decoupage is given by the fibered product of the puff groupoids $\gr V_{\mathcal{F}i}$ for each face. The groupoid $\Gamma_\phi(X)$ of a manifold with fibred corners X is then defined as the restriction of the groupoid of a decoupage in which X is embedded. One recovers the groupoid of $b$-calculus by Monthubert in the special case of a trivial fibration. The section ends with the proof of amenability and longitudinal smoothness of the groupoid $\Gamma_\phi(X)$.

\smallskip

\subsubsection{Definition}

\begin{definition}
Let $\eq_\phi=(M,(V_i,\phi_i)_{i\in I})$ be a fibred decoupage. The groupoid of $\eq_\phi$ is the fibered product of the puff groupoids $\gr V_{i\phi}$ through the maps $s\oplus r:\gr V_{i\phi}\rightarrow M^2$.
\end{definition}

\begin{definition}
Let X be a manifold with fibred corners and $\eq_\phi$ a fibred decoupage associated with X.
The groupoid of X is the restriction to X of the groupoid $\gr \eq_\phi$, 
$\Gamma_\phi(X)=(\gr \eq_\phi)_X^X$.
\end{definition}

\begin{remark}
Each $V_i$ splits M in two parts according to definition \ref{equar}.
\end{remark}
\begin{remark}
In the case of a manifold X with a connected fibred boundary $\phi:\partial X\rightarrow Y$, an expression of the groupoid of X as a set is:
$$\Gamma_\phi(X)=(\partial X\times_Y \partial X\times_YTY\times\Rr_+^*)\sqcup X^\circ\times X^\circ.$$
\end{remark}

Monthubert-Pierrot \cite{MP} and Nistor-Weinstein-Xu \cite{NWX} showed how to define a pseudodifferential calculus $\Psi_c^{\infty}(\gr)$ for any smooth differentiable longitudinally smooth groupoid $\gr$. 
We now prove that $\gr \eq_\phi$ is longitudinally smooth, which directly implies the longitudinal smoothness of 
$\Gamma_\phi(X)$. We can therefore associate to any manifold with fibred corners X a compact support calculus $\Psi_c^{\infty}(\Gamma_\phi(X))$. The properties of this calculus and the link with Melrose's $\phi$-calculus will be studied in the next section.

\begin{proposition}
The groupoid of a fibred decoupage is a Lie groupoid.
\end{proposition}

\begin{demo}
Denote $\varphi_i=s\oplus r:\gr V_{i\phi}\rightarrow M^2$.
Let $\gamma=(\gamma_i)_{i\in I}$ be an element of the fibered product $\gr \eq_\phi$ of the $\gr V_{i\phi}$.
Then there exists a subset J of I such that 
$$(\gamma_i)_{i\in I} \in \prod_{i\in J}\gr V_{i\phi|V_i}\prod_{i\not\in J}\gr V_{i\phi|(M\setminus V_i)}.$$
It is thus enough to show the transversality of the morphisms $\varphi_i':\gr V_{i\phi|V_i}\rightarrow M^2, i\in J$, which are the canonical projections, and $\phi_i':\gr V_{i\phi|(M\setminus V_i)}\rightarrow M^2$ which are the inclusions. The normal bundles of the latter being trivial, it suffices to consider that $i\in J$.
But the canonical projections factorize through the submersion $\psi_i: \gr V_{i\phi|V_i}\rightarrow 
V_i\times_{Y_i}V_i$: if $i_i$ is the inclusion of $V_i\times_{Y_i}V_i$ in $M^2$, then $\varphi_i'=i_i\circ\psi_i$. The transversality of the morphisms $i_i$ is implied by the definition of a fibred decoupage and the surjectivity of the differentials $d\psi_i$ then shows that the morphisms  $\varphi_i',i\in J$, and consequently $\varphi_i$, are transverse.

Thus $\gr \eq_\phi$ is a submanifold of $\prod_{i\in I}\gr V_{i\phi}$ naturally endowed with the composition law  induced by the inclusion, it is a Lie groupoid.
\end{demo}

\subsubsection{Longitudinal smoothness of the groupoid}
\begin{proposition}
\label{long_eclat}
The groupoid of a fibred decoupage is longitudinally smooth.
\end{proposition}
\begin{demo}
Let $x\in M$ and $F$ be the open face including $x$. Let $J\subset I$ such that $F$ is the interior of $\bigcap_{i\in J}V_j$. The fibre $\gr \eq_{\phi}^x$ of $\gr \eq_{\phi}$ in $x$ is composed as the fibered product :
$$\prod_{i\in J}\gr V_{i\phi|V_i}^x\prod_{i\not\in J}\gr V_{i\phi|(M\setminus V_i)}^x.$$
$\gr V_{i\phi|(M\setminus V_i)}^x=\emptyset$ when $i\not\in J$, as $x$ is an element of $\bigcap_{i\in J}V_j$.
Besides when $i\in J$, $\gr V_{i\phi|V_i}^x$ is the smooth vector bundle  $\mathscr{N}_{V_i\times_{Y_i}\{x\}}^{V_i\times V_i\times \Rr_+^*}$. 
The transversality of the maps  $s\oplus r:\gr V_{i\phi}\rightarrow M^2$ restricted to $\gr V_{i\phi}^x$ then implies the smoothness of the fibered product $\prod_{i\in J}\gr V_{i\phi|V_i}^x$ and therefore that of the fibre  $\gr \eq_{\phi}^x$.
\end{demo}

\subsubsection{Amenability of the groupoid}
\begin{proposition}
The groupoid of a fibred decoupage is amenable.
\end{proposition}
\begin{demo}
Recall (see \cite{Ren}) that if $\gr$ is a groupoid and $U$ an open set of $\gr$, then $\gr$ is amenable if and only if $\gr_U$ and $\gr_{X\setminus U}$ are amenable.
Also a groupoid such that $\varphi=s\oplus r:\gr\rightarrow (\go)^2$ is surjective and open is amenable if and only if its isotropy subgroups are amenable.

Since then $\gr V_\phi$ is amenable: the isotropy of $\gr V_{\phi|M\setminus V}=(M\setminus V)^2$ is trivial  and $\gr V_{\phi |V}$ has isotropy groups isomorphic to $T_yY\times\Rr$, which is abelian and thus amenable.

$\gr \eq_\phi$ is therefore also amenable as the fibered product of amenable groupoids.
\end{demo}

\section{Pseudodifferential calculus on manifolds with fibred corners}
We briefly recall the definition of $\phi$-calculus, a particular pseudodifferential calculus introduced by Mazzeo and Melrose on manifolds with fibred boundaries \cite{MM}. It is a generalization of the $b$-calculus, built from a geometric desingularization of the manifold with corner $X^2$ and the fibrations of its boundaries.

We then associate to any manifold with fibred corners X a pseudo\-differential calculus with compact support $\Psi_c^{\infty}(\Gamma_\phi(X))$ from the groupoid $\Gamma_\phi(X)$ built in sect\-ion 2 and we show that $\phi$-calculus identifies with the pseudodifferential calculus on this groupoid in the case of manifolds with fibred boundary.

Finally we introduce a Schwartz type algebra $\mathcal{S}_\psi(\Gamma_\phi(X))$ which allows to define an extended calculus $\Psi^{\infty}(\Gamma_\phi(X)) = \Psi_c^{\infty}(\Gamma_\phi(X)) + \mathcal{S}_\psi(\Gamma_\phi(X))$. We show that $\mathcal{S}_\psi(\Gamma_\phi(X))$ is stable under holomorphic functionnal calculus and that the operators of extended $\phi$-calculus are elements of this algebra.

\subsection{$b$-stretched product and $b$-calculus}
Following the work of H\"{o}rmander who showed that pseudodifferential calculus on a manifold $M$ is given by Schwartz kernel on $M^2$, Melrose defined $b$-calculus on manifolds with boundary using Schwartz kernels on a space playing the role of $M^2$, the ``$b$-stretched product".  

Let $X$ be a compact smooth manifold with boundary. The $b$-stretched product of $X$ (\cite{Mel}) is a compact smooth manifold with corners defined by the gluing 
$$X^2_b=X^2\setminus (\bX)^2\cup S_+N,$$
 where $S_+N=(N_{(\bX)^2}^{X^2}-\{0\})/\Rr_+^*$ is a fibred space which can be made trivial under the form $(\bX)^2\times [-1,1]$ by choosing a definition function $\rho$ of $\bX$.
The embedding of $X^2\setminus (\bX)^2$ in $X^2\times [-1,1]$ by
$$ (x,y)\rightarrow (x,y,\frac{\rho(x)-\rho(y)}{\rho(x)+\rho(y)})$$
allows the identification of $X^2_b$ with a submanifold of $X^2\times [-1,1]$.
Another possible embedding in $X\times X\times\Rr$ is given by :
\label{x2b}
$$X^2_b=\left \{ (x,y,t)\in X\times X\times \Rr,\; (1-t)\rho(x)=(1+t)\rho(y) \right \}.$$
The $b$-stretched product has three boundary components,
$lb = \{(x,y,-1)\in {\bX\times X\times\Rr}\}\simeq \bX\times X$, $rb = \{(x,y,1)\in X\times\bX\times\Rr\}\simeq X\times\bX$ and $S_+N$.

\begin{figure}[hbtp]
\begin{minipage}[b]{0.3\linewidth}
\centering
	\begin{pspicture}(0.5,-0.5)(3,4)
    \psset{linewidth=1pt}
    \psline{}(0.7071,0.7071)(3,3)	\uput[45](1.5,1){\LARGE\textbf{$\Delta_b$}}
	\psline{}(1,0)(2,0)	\uput[180](0,1.5){\Large\textbf{$lb$}}
	\psline{}(0,1)(0,2) \uput[-90](1.5,0){\Large\textbf{$rb$}}
	\psline[linestyle=dashed]{}(2,0)(3,0)
	\psline[linestyle=dashed]{}(0,2)(0,3)
    \psarc[linecolor=blue](0,0){1}{0}{90}
	\uput[-135](0.90,0.75){\Large\textbf{\textcolor{blue}{$S_+N$}}}
    \uput[180](1,2.5){\Large\textbf{$X^2_b$}}
	\end{pspicture}
\end{minipage}

\label{coin}
\caption{The $b$-stretched product of $R_+$.}
\end{figure}
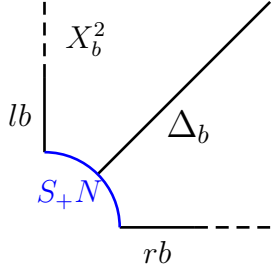

The diagonal $\Delta_b=\{(x,x,0)\in X\times X\times\Rr\}$ only intersects the boundary component $S_+N$ and this intersection is transverse, which allows a direct definition of operators with Schwartz kernels on $X^2_b$, the operators of $b$-calculus (see \cite{Mel}). The kernel $K$ of a $b$-pseudodifferential operator is an element of  $I^\infty(X^2_b,\Delta_b,\Omega^{\frac{1}{2}})$ with Taylor series development vanishing on $lb\cup rb$ ($K$ is a $\ci$ function in the neighbourhood of $lb\cup rb$ as the intersection of $lb\cup rb$ with $\Delta_b$ is empty).
The situation for $X=\Rr_+$ is illustrated on figure \ref{coin}.

\subsection{$\phi$-stretched product and $\phi$-calculus}
When the boundary of $X$ is the total space of a fibration $\phi:\bX\rightarrow Y$, the construction of $\phi$-calculus not only needs transverse intersection of the diagonal $\Delta_b$ with the $S_+N$ component but also with the submanifold $\Phi$ of $S_+N$ defined by
$$\Phi=\{(x,y,0)\in S_+N\simeq (\bX)^2\times [-1,1],\; \phi(x)=\phi(y)\}.$$

This condition is satisfied by introducing the $\phi$-stretched product \cite{MM}, which is obtained as a gluing 
$$X^2_\phi = X^2_b\setminus \Phi \cup S_+N_\phi$$
where $S_+N_\phi = (N_\Phi^{X^2_b}-\{0\})/\Rr_+^*$  is a fibred space which can be made locally trivial under the choice of a definition function $\rho$ of $\bX$ and a local trivialization $\tilde{\phi}$ of $\phi$ over $\bX$.



Boundary components of the $\phi$-stretched product are $lb$, $rb$, $S_+N_\phi$ and $sb$, the closure of $S_+N\setminus\Phi$ in $X^2_\phi$. When ${\dim Y = q}$, $X^2_\phi\setminus lb\cup rb\cup sb$ is shown to be diffeomorphic to the subset of elements $(x,x',S,Y)\in X\times X\times \Rr\times \Rr^q$ such that:
$$  \rho(x')=\rho(x)(1 + S\rho(x)) \; \text{ and } \; \tilde{\phi}(x) = \tilde{\phi}(x')-\rho(x)Y $$

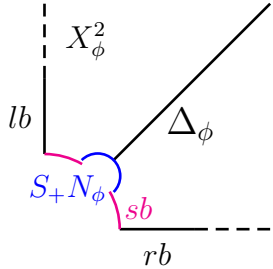
\begin{figure}[hbtp]
\label{pic_coins}
\begin{minipage}[b]{0.3\linewidth}
\centering
	\begin{pspicture}(0.5,-0.5)(3,4)
    \psset{linewidth=1pt}
    \psline{}(0.9192,0.9192)(3,3)	\uput[45](1.5,1){\LARGE\textbf{$\Delta_\phi$}}
	\psline{}(1,0)(2,0)	\uput[180](0,1.5){\Large\textbf{$lb$}}
	\psline{}(0,1)(0,2) \uput[-90](1.5,0){\Large\textbf{$rb$}}
	\psline[linestyle=dashed]{}(2,0)(3,0)
	\psline[linestyle=dashed]{}(0,2)(0,3)
    \psarc[linecolor=magenta](0,0){1}{0}{30}
    \psarc[linecolor=magenta](0,0){1}{60}{90}
    \psarc[linecolor=blue](0.7071,0.7071){0.3}{-45}{135}
    \uput[-90](1.25,0.6){\Large\textbf{\textcolor{magenta}{$sb$}}}
    \uput[-135](0.95,0.85){\Large\textbf{\textcolor{blue}{$S_+N_\phi$}}}
    \uput[180](1,2.5){\Large\textbf{$X^2_\phi$}}
	\end{pspicture}
\end{minipage}

\caption{The $\phi$-stretched product of $R_+$ with the trivial fibration.}
\end{figure}

The diagonal $\Delta_\phi=\{(x,x)\in X^2_\phi\}$ only intersects the boundary component $S_+N_\phi$ and this intersection is transverse  which allows a direct definition of operators with Schwartz kernels on $X^2_\phi$, the operators of $\phi$-calculus (see \cite{MM}).
The definition is similar to the one of previous paragraph. The kernel $K$ of a $\phi$-pseudodifferential operator is an element of $I^\infty(X^2_\phi,\Delta_\phi,\Omega^{\frac{1}{2}})$  with Taylor series development vanishing on $lb\cup rb\cup sb$ ($K$ is a $\ci$ function in the neighbourhood of $lb\cup rb\cup sb$ as the intersection of $lb\cup rb\cup sb$ with $\Delta_\phi$ is empty).
The situation for $X=\Rr_+$ is illustrated on figure \ref{pic_coins}.

\subsection{$\phi$-calculus and the groupoid of manifolds with fibred boundary}

From the previous discussions it appears that two geometric objects model calculus on manifolds with fibred boundary : on the one hand $\phi$-stretched product $X^2_\phi$ and the associated $\phi$-calculus, on the other hand the groupoid $\Gamma_\phi(X)=(\gr \eq_\phi)_X^X$ and its pseudodifferential calculus $\Psi_c^{\infty}(\Gamma_\phi(X))$.

We show here that the groupoid $\Gamma_\phi(X)$ of a manifold with fibred boundary can be embedded in the $\phi$-stretched product $X^2_\phi$ as the open submanifold $X^2_\phi\setminus lb\cup rb\cup sb$  and that through this identification $\phi$-calculus coincides with the pseudodifferential calculus $\Psi_c^\infty(\Gamma_\phi(X))$.

\subsubsection{Notations}
Let X be a manifold with fibred boundary $\phi:\bX\rightarrow Y$,  $p=\dim{\bX}-\dim{Y}$ and $q=\dim{Y}$.
Denote $\bX_j$ the connected components of $\bX$ with induced fibrations $\phi_j:\bX_j\rightarrow Y_j$, and $\eq_\phi=(M,(\bX_j,\phi_j)_{j\in J})$ the fibred decoupage which positive part defines X and its fibration.
As the connected components of $\bX_j$ are disjoint, the puff $\gr \eq_\phi$ of $\eq_\phi$ coincides with the union $\cup_{j\in J}\gr \bX_{j\phi}$ of the puff groupoids of $\bX_j$ in $M$. It is given below an explicit $\ci$ atlas of $\gr \eq_\phi$ on which the embedding is defined.

Let $\{(\mathcal{U}_{\alpha}=\mathcal{W}_{\alpha}\times\Rr_+^*,\varphi_{\alpha}=\psi_{\alpha}\times\log)\}$ be a $\ci$ atlas of $\bigcup_{j\in J}\bX_j\times\bX_j\times\Rr_+^*$ which charts are slices of each  $\mathcal{H}ol(\phi_j)\times\Rr_+^*=\bX_j\times_{Y_j}\bX_j\times\Rr_+^*$. For a given  $\mathcal{U}=\mathcal{U}_\alpha$, $U = \varphi_{\alpha}(\mathcal{U}_{\alpha})$ is an open set of $\Rr^{2p+q}\times\Rr^q\times\Rr$, a slice of $V={U\cap (\Rr^{2p+q}\times\{0\}\times\{0\})}$ and the composition 
$$ \mathscr{D}_\mathcal{V_\alpha}^\mathcal{U_\alpha} \stackrel{\tilde{\varphi}_\alpha}{\longrightarrow} \mathscr{D}_{V}^{U} \stackrel{\Theta^{-1}}{\longrightarrow}\Omega_{V}^{U} $$ 

\noindent is a diffeomorphism of $\mathscr{D}_\mathcal{V}^\mathcal{U}$ on the open set $\Omega_{V}^{U}$ of $\Rr^{2p+q}\times\Rr^q\times\Rr\times\Rr$ defined by (see also section \ref{cone_diff}):
$$\Omega_{V}^{U} = \{ (x,\xi,\lambda,t) \in \Rr^{2p+q}\times\Rr^q\times\Rr\times\Rr \; , \; (x,t\xi,t\cdot e^\lambda)\in U \}.$$

\subsubsection{Embedding of $\Gamma_\phi(X)$ in $X^2_\phi$}
Let $F:U\rightarrow \Rr^{2p+q}\times\Rr^q\times\Rr$ be the $\ci$ map defined by $(x,\xi,\lambda)\mapsto (x,\xi,e^\lambda-1)$. Let denote $U'=F(U)$ and $V'=F(V)$. The relation $F(x,0,0)=(x,0,0)$ implies the map $\tilde{F}:\Omega_{V}^{U}\rightarrow \Omega_{V'}^{U'}$ (see \ref{cone_diff}) to be a diffeomorphism.

Then consider the map $\Theta_\alpha:\mathscr{D}_\mathcal{V_\alpha}^\mathcal{U_\alpha}\rightarrow \Rr^{2p+q}\times\Rr^q\times\Rr$ defined by the composition
$$ \mathscr{D}_\mathcal{V_\alpha}^\mathcal{U_\alpha} \stackrel{\Theta^{-1}\circ\tilde{\varphi}_\alpha}{\longrightarrow} \Omega_{V_\alpha}^{U_\alpha}  \stackrel{\tilde{F}}{\longrightarrow} \Omega_{V_\alpha'}^{U_\alpha'}$$
and its projection $\pi\Theta_\alpha:\mathscr{D}_\mathcal{V_\alpha}^\mathcal{U_\alpha}\rightarrow \Rr^q\times\Rr$ over the last two factors.
If $s$ and $r$ denote the source and target maps of $\Gamma_\phi(X)$ we define ${i_s:\Gamma_\phi(X)\rightarrow X\times X\times\Rr^q\times\Rr}$ for $\gamma\in i_\D(\mathscr{D}_\mathcal{V_\alpha}^\mathcal{U_\alpha})$ by :
$$i_s(\gamma)= \left ( s(\gamma),r(\gamma),\pi\Theta_\alpha\circ i_\D^{-1}(\gamma) \right )$$

\begin{proposition}
\label{plongement}
The map $i_s$ is an embedding of $\Gamma_\phi(X)$ in $X^2_\phi$.
\end{proposition}
\begin{demo}
Let denote $\tilde{\phi}:\mathcal{W}_{\alpha}\rightarrow \Rr^q$ the projection of $\psi_\alpha$ on the $\Rq$ factor of $\psi_\alpha(\mathcal{W}_{\alpha}) \subset \Rr^{2p+q}\times\Rr^q$.
The expression \ref{hol_diff} of $i_\D^{-1}$ for an element $(x,y)$ of $X^\circ\times X^\circ$ shows that  $$\tilde{\varphi}_\alpha\circ i_\D^{-1}(x,y)= ( \psi_\alpha(\pi(x),\pi(y)),\log{\frac{\rho(y)}{\rho(x)}},\rho(x) ).$$

The relation
$$F( \psi_\alpha(\pi(x),\pi(y)),\log{\frac{\rho(y)}{\rho(x)}} ) = \left ( \psi_\alpha(\pi(x),\pi(y)),\frac{\rho(y)-\rho(x)}{\rho(x)} \right )$$ then implies 
$$\pi\Theta_\alpha\circ i_\D^{-1}(x,y) = \left ( \frac{\tilde{\phi}(\pi(x),\pi(y))}{\rho(x)} , \frac{\rho(y)-\rho(x)}{\rho^2(x)} \right ).$$

Now $X^2_\phi\setminus lb\cup rb\cup sb$ is the submanifold of $X\times X\times\Rr^q\times\Rr$ with elements $(x,y,Y,S)$ such that $  \rho(y)=\rho(x)(1 + S\rho(x)) \; \text{ and } \; \tilde{\phi}(x) = \tilde{\phi}(y)-\rho(x)Y .$
The element $i_s(x,y)=(x,y,Y,S)$ where $Y= \frac{\tilde{\phi}(\pi(x),\pi(y))}{\rho(x)}$ and $S=\frac{\rho(y)-\rho(x)}{\rho^2(x)}$ thus belongs to $X^2_\phi\setminus lb\cup rb\cup sb$ as :
$$ \rho(x)(1 + S\rho(x)) = \rho(x)(1 + \frac{\rho(y)-\rho(x)}{\rho(x)} ) =  \rho(y)$$
and $$ \tilde{\phi}(y)-\rho(x)Y = \tilde{\phi}(y) - \tilde{\phi}(\pi(x),\pi(y)) = \tilde{\phi}(x).$$

The trivial holonomy of fibrations implies the injectivity of source and target maps of $\Gamma_\phi(X)$ and the injectivity of $i_s$ restricted to the interior of $\Gamma_\phi(X)$.
Moreover the fiber in 0 of a deformation $\mathscr{D}_\mathcal{V_\alpha}^\mathcal{U_\alpha}$ is diffeomorphic to the set $\Omega_{V0}^{U} = \{ (x,\xi,\lambda,0) \in \Rr^{2p+q}\times\Rr^q\times\Rr_+^*\times\{0\} \; , \; x\in V \}$, which proves the surjectivity of $\pi\Theta_\alpha\circ i_\D^{-1}$ when restricted to the boundary $\bigcup_{\alpha}\mathscr{D}_{\mathcal{V_\alpha}0}^\mathcal{U_\alpha}$ of $\Gamma_\phi(X)$.

Therefore $i_s(\Gamma_\phi(X)) \subset X^2$ is in bijection with $X^2_\phi\setminus lb\cup rb\cup sb$, the definition of $i_s$ as a composition of smooth maps allows to conclude that $i_s$ is an embedding of $\Gamma_\phi(X)$ in $X^2_\phi$, with image $X^2_\phi\setminus lb\cup rb\cup sb$.
\end{demo}

\subsubsection{Identification of $\Gamma_\phi(X)$ in $X^2_\phi$}
\label{identification}

\begin{proposition}
$\Gamma_\phi(X)$ is an open submanifold of $X^2_\phi$ and
$$ X^2_\phi\setminus\Gamma_\phi(X) = lb\cup rb\cup sb.$$
\end{proposition} 

\begin{demo}
The (non-disjoint) boundary components of $X^2_\phi$ are
$${\partial{X_\phi^2}=lb\cup rb\cup sb \cup S_+N_\phi.}$$
The relation $i_s(\Gamma_\phi(X)) = X^2_\phi\setminus lb\cup rb\cup sb$ established in the previous demonstration then proves that the closed submanifold $lb\cup rb\cup sb$ is the complementary in $X^2_\phi$ of the image of the embedding $i_s$.
\end{demo}

\subsection{Identification of $\Psi_c^\infty(\Gamma_\phi(X))$ with $\phi$-calculus}
Let denote $\Psi_c^\infty(\Gamma_\phi(X))$ the algebra of operators with compact support on $\Gamma_\phi(X)$ and $\Psi_{\phi,c}^\infty(X)$ the algebra of operators with compact support of $\phi$-calculus.
The submanifold $lb\cup rb\cup sb$ is disjoint from a neighbourhood of $\Delta_\phi$. Pseudodifferential calculus with compact support on $\Gamma_\phi(X)$ thus coincides with distributional sections on $X^2_\phi$ which vanish in a neighbourhood of $lb\cup rb\cup sb$. One recovers the definition of $\phi$-calculus :

\begin{theorem}
\label{phi_c}
The pseudodifferential calculus on $\Psi_c^\infty(\Gamma_\phi(X))$ coincides with the small $\phi$-calculus with compact support $\Psi_{\phi,c}^\infty(X)$ of Melrose.
\end{theorem}

\subsection{Extended pseudodifferential calculus}

We define an extended pseudodifferential calculus $\Psi^{\infty}(\Gamma_\phi(X))$ over the groupoid $\Gamma_\phi(X)$ of a manifold with fibred corners :   
$$\Psi^{\infty}(\Gamma_\phi(X)) = \Psi_c^{\infty}(\Gamma_\phi(X)) + \mathcal{S}_\psi(\Gamma_\phi(X)).$$

The algebra $\mathcal{S}_\psi(\Gamma_\phi(X))$ of functions with rapid decay built from a length function  $\psi$ with polynomial growth (lemma \ref{length_f}) is naturally stable under holomorphic functionnal calculus. It includes the regularizing operators of $\phi$-calculus (proposition \ref{phi_ext}) and the extended $\phi$-calcul satisfies the inclusion relation : $$\Psi_\phi^\infty(X)\subset\Psi^\infty(\Gamma_\phi(X)).$$

\smallskip

\subsubsection{Length functions with polynomial growth}

\begin{definition}[\cite{Mont}, 1.4]
Let $\gr$ be a Lie groupoid and $\mu$ a Haar system on $\gr$. A length function with polynomial growth is a $\ci$ function $\varphi:\gr\rightarrow \Rr_+$ such that :
\begin{itemize}
\item[$\cdot$] $\varphi$ is subadditive, i.e. $\varphi(\gamma_1\gamma_2)\leq \varphi(\gamma_1) + \varphi(\gamma_2)$,
\item[$\cdot$] $\forall \gamma\in G, \varphi (\gamma^{-1})=\varphi (\gamma)$,
\item[$\cdot$] $\varphi$ is proper,
\item[$\cdot$] $\exists c,N,\forall x\in G^{(0)},\forall r\in\Rr_+,\mu_x(\varphi^{-1}([0,r]))\leq c(r^N + 1)$.
\end{itemize}
\end{definition}

With such a function $\varphi$ one can define the space 
$$ \mathcal{S}^0=\{ f\in C_0(\gr,\Omega^{1/2}),\forall P\in\mathbb{C}[X],\sup_{\gamma\in \gr}|P(\varphi(\gamma))f(\gamma)|<\infty \} $$
where the half-densities bundle $\Omega^{\frac{1}{2}}$ is the line bundle over $\gr$ which fibre for $\gamma\in\gr$ is the vector space of applications $$\rho:\Lambda^kT_\gamma \gr^{r(\gamma)}\otimes\Lambda^kT_\gamma \gr_{s(\gamma)}\rightarrow \mathbb{C}$$ such that $\rho(\lambda\nu)=|\lambda|^{1/2}\rho(\nu),\forall\lambda\in\Rr$.

The space $\mathcal{S}_\varphi(G)$ of functions with rapid decay over $\gr$ is defined in (\cite{Mont}, 1.5) as the subspace of $\mathcal{S}^0$ of functions $f$ such that :
$$\forall l\in\mathbb{N},\forall (v_1,\dots,v_l)\in \ci(\mathcal{A}G)^l,\forall k\leq l,(v_1\dots v_k\cdot f\cdot v_{k+1}\dots v_l)\in \mathcal{S}^0 .$$

$\mathcal{S}_\varphi(G)$ is the ideal of regularizing operators of the extended pseudodifferential calculus  $I^\infty(G,G^{(0)};\Omega^{1/2}) = I_c^\infty(G,G^{(0)};\Omega^{1/2}) + \mathcal{S}_\varphi(G)$. It is also a subalgebra of $C^*(\gr)$ stable under holomorphic functionnal calculus (\cite{LMN2}, theorem 6).

\smallskip

\subsubsection{Length function for manifolds with fibred corners}
Let $X$ be a manifold with embedded and fibred corners defined by a fibred decoupage $\eq_\phi = (M,(V_i,\phi_i)_{i\in [1,N]})$. Each puff groupoid $\gr V_{i\phi}$ is diffeomorphic through the embedding $i_s$ of proposition \ref{plongement} to a closed submanifold of $M\times M\times\Rr^{q_i}\times\Rr$. The positive part of their fibered product over $M\times M$ is thus diffeomorphic (through a diffeomorphism $i$) to a closed submanifold of $X\times X\times \prod_{i\in [1,N]}\Rr^{q_i}\times\Rr^N$. Let $\pi: X\times X\times \prod_{i\in [1,N]}\Rr^{q_i}\times\Rr^N \rightarrow \prod_{i\in [1,N]}\Rr^{q_i}\times\Rr^N$ be the canonical projection, we define a length function on $\gr\eq_\phi$ by $\psi(\gamma)=\norm{\pi\circ i(\gamma)}.$

\begin{lemma}
\label{length_f}
$\psi:\gr\eq_\phi\rightarrow \Rr_+$ is a length function with polynomial growth.
\end{lemma}
\begin{demo}
On every chart $i_\D(\D_{\mathscr{V}_\alpha}^{\mathscr{U}_\alpha}) \subset \gr V_{i\phi}$ of proposition \ref{plongement}, 
$\psi_i=\norm{\pi\Theta\circ i_\D^{-1}(\gamma)}$ satisfies $\psi_i(\gamma)=\norm{\xi}+\norm{\lambda}$ where $(x,\xi,\lambda,t)=\Theta^{-1}\circ \tilde{\phi}\circ i_\D^{-1}(\gamma)$.
Each $\psi_i$ is therefore a groupoid morphism and $\psi=\sum_{i\in [1,N]}\psi_i$ is subadditive.

$\pi$ is proper as X is compact, besides $i(\Gamma_\phi(X))$ is closed in $X\times X\times \prod_{i\in [1,N]}\Rr^{q_i}\times(\Rr)^N$ so $\psi=\norm{\pi\circ i}$ is proper.

Finally for $x\in X$ the local expresion of the $\psi_i$ shows that there exists two constants $c_1,c_2$ such that for $r\geq 1$, $\mu_x(\psi^{-1}([0,r])) < c_1\vol(B_{\Rr^M}(0,r))$ where $M=N+\sum_iq_i$ and for $r<1$, $\mu_x(\psi^{-1}([0,r])) \leq \sup_{x\in X}\mu_x(\psi^{-1}([0,1]))=c_2.$
By denoting $c=\max(c_1,c_2)$ we get $\mu_x(\psi^{-1}([0,r]))\leq c(r^M+1)$. As $X$ is compact it can be covered with a finite number of neighbourhoods which give an inequality on $X$.

Hence $\psi$ has polyomial growth.
\end{demo}

\subsubsection{Definition of the extended pseudodifferential calculus}
We are now able to define an extended pseudodifferential calculus on manifolds with fibred corners:
\begin{definition}
Let X be a manifold with fibred corners and $\Gamma_\phi(X)$ the associated longitudinally smooth groupoid. The extended pseudodifferential calculus on $X$ is the algebra $\Psi^{\infty}(\Gamma_\phi(X))$ of pseudodifferential operators on $\Gamma_\phi(X)$ defined by:
 $$\Psi^{\infty}(\Gamma_\phi(X)) = \Psi_c^{\infty}(\Gamma_\phi(X)) + \mathcal{S}_\psi(\Gamma_\phi(X))$$ 
\end{definition}

\noindent where $\mathcal{S}_\psi(\Gamma_\phi(X))=\mathcal{S}_\psi((\gr \eq_\phi)_X^X)$ is the space of functions with rapid decay with respect to the length function $\psi$ of previous lemma. 

$\mathcal{S}_\psi(\Gamma_\phi(X))$ is the ideal of regularizing operators of the extended pseudodifferential calculus $\Psi^{\infty}(\Gamma_\phi(X))$. 

\begin{remark}
The definition of $\mathcal{S}_\psi(\Gamma_\phi(X))$ given here is based on the existence of a length function with polynomial growth. In more complex cases of foliated boundary no such functions exist in general, for instance in the case of nonzero entropy foliations. The problem is overcome in section 3.5 of \cite{Gui} with the definition of a Schwartz algebra $\mathcal{S}_c(\Gamma_\phi(X))$ valid for any regular foliation of the boundary. For the sake of clearness only the fibration case is developped in this article.
\end{remark}

\smallskip

\subsubsection{Identification of $\Psi^\infty(\Gamma_\phi(X))$ and $\Psi_\phi^\infty(X)$}
The regularizing operators of $\phi$-calculus are $\ci$ kernels on $X^2_\phi$ which Taylor series vanish on $lb\cup rb\cup sb$. Melrose shows in proposition 4 of \cite{MM} that a kernel $K(x,x')$ which vanishes in Taylor series with respect to the powers of $\rho(x),\tilde{\phi}(x)$ is a rapid decay kernel with respect to the variables $Y=\frac{\tilde{\phi}(x)-\tilde{\phi}(x')}{\rho(x)}$ and $S=\frac{\rho(x)-\rho(x')}{\rho(x)}$, that is with respect to $\pi\Theta\circ i_\D^{-1}(x,x')$. 	

Those kernels are thus \textsl{a fortiori} functions with rapid decay with respect to $\psi$, so that any regularizing operator of $\phi$-calculus defines by restriction to the complementary of $lb\cup rb\cup sb$ a function with rapid decay on $\Gamma_\phi(X)$ and $\Psi_{\phi}^{-\infty}(X)\subset \mathcal{S}_\psi(\Gamma_\phi(X))$.

We get as a corollary the following result:

\begin{theorem}
\label{phi_ext}
The extended pseudodifferential calculus $\Psi^\infty(\Gamma_\phi(X))=\Psi_c^{\infty}(\Gamma_\phi(X)) + \mathcal{S}_\psi(\Gamma_\phi(X))$ contains the operators of $\phi$-calculus :
$$\Psi_\phi^\infty(X)\subset\Psi^\infty(\Gamma_\phi(X))$$
\end{theorem}

\subsubsection{Identification of the boundary restriction with Melrose's normal operator}
Contrary to the case of a smooth manifold without corners, regularizing operators are not necessarly compact. The obstruction to compacity for $\phi$-calculus is described by introducing an indicial algebra $\mathit{I}(A,\xi,\lambda)\in \Psi_{\text{sus}(\phi)}(\partial X)$ and a normal operator $N_\phi:\Psi_\phi^\infty(X)\rightarrow \Psi_{\text{sus}(\phi)}(\partial X)$ with values in a suspended algebra of kernels over $\partial X$ which are invariant under translation and with rapid decay at infinity (see \cite{Mel}, \cite{MM}).

In the case of manifolds with fibred corners, the indicial family $\mathit{I}(A,\xi,\lambda)$ is defined by Melrose (\cite{MM}, proposition 5) as the restriction of the kernel to $ff(X_\phi^2)$. The normal operator $N_\phi$ of Melrose thus coincides with the operator $\partial$ of restriction to the boundary of the groupoid $\Gamma_\phi(X)$.
The translation invariance of the indicial family seen as an element of $\partial\mathcal{S}_\psi(\Gamma_\phi(X))$ is nothing but the invariance under $TY\times\Rr$ of a family of operators on $\mathcal{S}_\psi(\partial X\times_Y \partial X\times_Y TY\times\Rr)\simeq \mathcal{S}(TY\times\Rr,\cc(\partial X\times_Y \partial X))$.

\smallskip

\subsection{Total ellipticity and Fredholm index}
Let $\overline{\Psi}_0(\Gamma_\phi(X))$ be the norm closure of $\Psi_c^0(\Gamma_\phi(X))$ in the multiplier algebra of $C_r^*(\Gamma_\phi(X))$.
Let recall that the symbol map $\sigma$ induces the exact sequence:
$$ 0\rightarrow C_r^*(\Gamma_\phi(X))\longrightarrow \overline{\Psi}_0(\Gamma_\phi(X)) \stackrel{\sigma}{\longrightarrow} C_0(S^*\Gamma_\phi(X)) \rightarrow 0 $$

\noindent and that the analytic index $\ind_a:K^0(A^*\Gamma_\phi(X))\rightarrow K_0(C_r^*(\Gamma_\phi(X)))$ takes values in the $K$-theory of the $C^*$-algebra of the groupoid $\Gamma_\phi(X)$. The analytic index is thus not a Fredholm index in general, its values are not elements of $\mathbb{Z}$. To get a Fredholm operator it is necessary to introduce additionnal conditions, classical ellipticity (the property of a symbol to be invertible) not being sufficient. 

The condition of total ellipticity is derived very naturally from the groupoid approach. Indeed $\partial\Gamma_\phi = (\Gamma_\phi(X))_{\bX}^{\bX}$ is a closed saturated set of $\Gamma_\phi(X)$ and the following diagram is commutative:
\[
\xymatrix{               
0 \ar[r] & C_r^*(X^\circ\times X^\circ) \ar[d] & & & \\
0 \ar[r] & C_r^*(\Gamma_\phi(X)) \ar[d]\ar[r] & \overline{\Psi}_0(\Gamma_\phi(X)) \ar[d]^{\partial}\ar[r]^{\sigma} & C_0(S^*\Gamma_\phi(X)) \ar[d]\ar[r] & 0 \\
0 \ar[r] & C_r^*(\partial\Gamma_\phi) \ar[d]\ar[r] & \overline{\Psi}_0(\partial\Gamma_\phi) \ar[r] & C_0(S^*\partial \Gamma_\phi) \ar[d]\ar[r] & 0 \\
 & 0 & & 0 &
}
\] 

The relation $C_r^*(X^\circ\times X^\circ)=\ker(\sigma)\cap\ker(\partial)=\ker(\sigma\oplus\partial)$ shows that  $\sigma\oplus\partial$ factorizes through the map $\sigma_{\tau}$:
$$\overline{\Psi}_0(\Gamma_\phi(X)) \stackrel{\sigma_{\tau}}{\longrightarrow} C_0(S^*\Gamma_\phi(X)) \times_{\partial X} \overline{\Psi}_0(\partial\Gamma_\phi)$$
where $C_0(S^*\Gamma_\phi(X)) \times_{\partial X} \overline{\Psi}_0(\partial\Gamma_\phi)$ denotes the fibered product of $C_0(S^*\Gamma_\phi(X))$ et de $\overline{\Psi}_0(\partial\Gamma_\phi)$ over $C_0(S^*\partial \Gamma_\phi)$.

\begin{definition}[\cite{Gui}, 4.6]
An operator $P\in \overline{\Psi}_0(\Gamma_\phi(X))$ will be called \textsl{totally elliptic} when the element $\sigma_{\tau}(P)$ is inversible.
\end{definition}

The map $\sigma_{\tau}$ induces from the previous diagram the exact sequence :
$$ 0\rightarrow C^*(X^\circ\times X^\circ)\simeq \mathscr{K}\rightarrow \overline{\Psi}_0(\Gamma_\phi(X)) \stackrel{\sigma_{\tau}}{\rightarrow} C_0(S^*\Gamma_\phi(X)) \times_{\partial X} \overline{\Psi}_0(\partial\Gamma_\phi) $$

Therefore an element in $\overline{\Psi}_0(\Gamma_\phi(X))$ is invertible modulo compact operators if an only if its image in $C_0(S^*\Gamma_\phi(X)) \times_{\partial X} \overline{\Psi}_0(\partial\Gamma_\phi)$ is invertible. This last condition is equivalent to the double condition of ellipticity for the operator and inversibility for its bounday restriction.

In particular we get the property stated in the first part of \cite{MRo} : an operator of $\phi$-calculus is totally elliptic when its symbol and the normal operator $N_\phi=\partial$ are jointly invertible, and the operator is Fredholm if and only if it is totally elliptic.

\section{The holonomy groupoid of manifolds with fibred corners}
In section 2 we constructed the puff groupoid $\gr V_\mathcal{F}$ of a foliated submanifold $(V,\mathcal{F})$ of codimension 1 in M. The module $\F$ of vector fields tangent to $\mathcal{F}$ on V defines a singular foliation on M. We prove here that $\gr V_\mathcal{F}$ is the holonomy groupoid $\HolM$ of the singular foliation $(M,\F)$.

For that purpose we show the groupoid $\gr V_\mathcal{F}$ integrates the Lie algebroid $\A\F$ of $\F$ and satisfies the minimality condition defining the groupoid $\HolM$. Notations from section \ref{notations} are reused.

\subsection{The singular foliation of a fibred boundary}
Consider the Lie algebroid $\A\F=TM$ with anchor map $p$ defined over each $N_i$ by:
$$(f_i)_*(p((f_i)_*^{-1}(v_1,v_2,t,\lambda)))=(v_1+tv_2,t,t\lambda)$$
when $(v_1,v_2,(t,\lambda))\in TF_{i} \oplus \mathscr{N}_i \times T_t(\Rr)$, and by the identity over  $T(M\setminus N)$.

$\A\F$ defines a foliation $(M,\F)$ with this map (remark \ref{algfol}). Let recall that a foliation is \textsl{almost regular} when it is defined by a Lie algebroid with anchor map injective over a dense open set.

\begin{proposition}
The foliation $(M,\F)$ is almost regular.
\end{proposition}
\begin{demo}
It is sufficient to show that the anchor map of $\A\F$ is injective over the dense open set $M \setminus V$, which follows immediatly from the injectivity of the map $(t,\lambda)\mapsto (t,t\lambda)$ for $t\neq 0$.
\end{demo}

A result by Claire Debord (\cite{Deb},Thm 4.3) valid for any almost regular foliation then ensures the existence of a Lie groupoid $\HolM$ which Lie algebroid is $\A\F$.
It is a groupoid with units M which integrates the foliation defined by $\A\F$.

The immediate surjectivity of the map $(t,\lambda)\mapsto (t,t\lambda)$ for $t\neq 0$, and thus that of the anchor map of  $\A\F$ over $M \setminus V$ implies the leaves of this foliation are $M \setminus V$ and the leaves of $\mathcal{F}$. The groupoid $\grVF$ integrates this space of leaves and one can ask if it is the only one. The following property answers that question. 

A Lie algebroid $\mathscr{A}$ is said \textsl{extremal} for the foliation $\mathcal{F}$ if for any other Lie algebroid $\mathscr{A}'$ defining $\mathcal{F}$, there exists a Lie algebroid morphism from $\mathscr{A}'$ to $\mathscr{A}$. When the foliation $\mathcal{F}$ is almost regular, the extremality of $\mathscr{A}$ is equivalent (\cite{Deb}) to the existence of a unique ``minimal" groupoid G integrating $\mathscr{A}$, in the sense that for any other Lie groupoid $H$ defining $\mathcal{F}$ there exists a differentiable morphism of groupoids from $H$ to $G$. In that sense $\grVF$ can be defined as the ``smallest" Lie groupoid defining  $\mathcal{F}$. Indeed:

\begin{proposition}
$\A\F$ is extremal for $(M,\F)$.
\end{proposition}
\begin{demo}
It is sufficient to show that the anchor map $p$ of $\A\F$ induces an isomorphism between $\Gamma(\A\F)$, the vector space of local sections of $\A\F$ and $\Gamma(T\F)$, the vector space of local vector fields tangent to $\F$.
The condition over $M\setminus V$ immediately results from the bijectivity of the anchor map $p_{|M\setminus V}$.
The condition over $V$ is obtained from the Taylor series development of a vector field tangent on V to the leaves of the foliation $\mathcal{F}$. If $x_i,y_i,t$ denote local coordinates associated to the directions respectively defined by $TF_{i}$, $\mathscr{N}_i$ and $T_t(\Rr)$, a vector field  $\xi\in\F$ will be locally generated by the free family $\langle t\frac{\partial}{\partial t},\frac{\partial}{\partial x_i},t\frac{\partial}{\partial y_i} \rangle$, which ensures the announced isomorphism.
\end{demo}

\subsection{The groupoid of a fibred decoupage as a holonomy groupoid}

\begin{definition}
We call holonomy groupoid of $(M,\F)$ the minimal Lie groupoid  $\mathcal{H}ol(M,\F)$ with Lie algebroid $\A\F$. The previous proposition shows that $\mathcal{H}ol(M,\F)$ is generated from the atlas of bi-submersions near the identity of vector fields tangent on V to the leaves of $\mathcal{F}$ (see \cite{ASk}).
\end{definition}

\begin{lemma}
\label{semi_regular}
$\grVF$ is a semi-regular s-connected Lie groupoid.
\end{lemma}

\begin{demo}
The s-connected property of $\grVF$ is immediate by definition.
The groupoid $\D_\varphi$ is endowed with the differential structure of a normal cone deformation. The open sets of  $G_V\times\Rr_+^*\times \Rr^*$ and the sets of the form $\Theta(\Omega_V^U)\exp{Z}$, $\Omega_V^U = \{ (x,x^\prime,\xi,t) \in \Rr^{2p}\times\Rq\times\Rr\times\Rr \; , \; (x,x^\prime,t\xi)\in U \}$ form a regular open sets basis for the topology of $\D_\varphi$. Thus $\D_\varphi$ is semi-regular and so is $\grVF$ according to proposition \ref{psi_G}.
\end{demo}

\begin{proposition}[\cite{Deb}]
\label{quasi_graphoid}
Let H be a semi-regular s-connected Lie groupoid. Then H is a quasi-graphoid if and only if the set $\{ x\in H^{(0)} \;,\; H_x^x=\{x\} \}$ of units with trivial isotropy group is dense in $H^{(0)}$.
\end{proposition}

\begin{corollary}
$\grVF$ is a quasi-graphoid.
\end{corollary}

\begin{demo}
It is immediate to check that the isotropy of $\grVF$ is trivial over $M^\circ$. $\gr V_\mathcal{F}$ being s-connected and semi-regular (lemma \ref{semi_regular}), it is a quasi-graphoid from proposition \ref{quasi_graphoid}.
\end{demo}

\begin{proposition}
\label{gr_holonomy}
$\grVF$ is the holonomy groupoid $\HolM$ of $(M,\F)$.
\end{proposition}

\begin{demo}
Let $\mathcal{A}$ and $\mathcal{B}$ be the atlas introduced in \ref{hol_diff} to describe the differential structure of $\gr V_\mathcal{F}$. $\mathcal{A} = \{(\Omega_i^j,\phi_{ij})_{i,j\in I^2}\}$ with $\phi_{ij}=\varphi_{ij}\circ i_\D^{-1}$ and $\Omega_i^j=i_\D(\D_i^j)$.
$\mathcal{B} = \{(\Omega_\beta,\varphi_\beta)\}$ is an atlas of the manifold $M^{\circ}_+\times M^{\circ}_+ \bigsqcup M^{\circ}_-\times M^{\circ}_-$.

Over $\Omega_\mathcal{B}=\bigcup\Omega_\beta$ source and target maps of $\gr V_\mathcal{F}$ are defined by:
\[ \left \{ \begin{array}{ccc} s(x,y) &=& x \\
                               r(x,y) &=& y
             \end{array} \right .
\]             
The algebroid $\A\grVF$ restricted to $M\setminus N\subset\Omega_\mathcal{B}$ is thus $\A_{M\setminus N}=\bigcup\ker(ds_{|\{(x,x),x\in M\setminus V\}})=T(M\setminus V)$ endowed with the identity anchor map $p:\A_{M\setminus N}\circ\rightarrow T(M\setminus V)$ as $dr_{|T(M\setminus V)}=\id$. So $\A\gr V_{\mathcal{F}|M\setminus N}=\A\mathcal{F}_{M|M\setminus N}$.

Besides the groupoid immersion $\varphi:G_1\rightarrow V\times V$ induces an injective Lie algebroid morphism 
$\varphi^{\mathscr{A}}_*:\mathscr{A}_1\rightarrow TV$, with $\mathscr{A}_1=TF$. The connection $\omega$ allows the identification of $TV$ with ${\A}_1\oplus\mathcal{N}$ where $\mathcal{N}=TV/{\A}_1$ and the anchor of the algebroid ${\A}_\varphi$ from the normal cone groupoid $D_\varphi$ is:
\[ \begin{array}{cccc}p_D:& \A_\varphi=\A_1\oplus\mathcal{N}\times\Rr & \rightarrow & T(V\times\Rr)= TV\times T\Rr \\
                                    & (v_1,v_2,t) & \mapsto & (p_2(v_1+tv_2),(t,0))
\end{array}
\]
The composition law on $\D_\varphi$ is given by the product of the law composition on $D_\varphi$ and the multiplicative action of $\Rr_+^*$ on $\Rr$. Therefore $\mathscr{A}\D_\varphi=\mathscr{A}_1\oplus\mathcal{N}\oplus T\Rr$ and the anchor map  $p_\varphi$ of $\mathscr{A}\D_\varphi$ is obtained as:
\[ \begin{array}{cccc}p_\varphi: & \mathscr{A}=\A_1\oplus\mathcal{N}\times T\Rr & \rightarrow & T(V\times\Rr)= TV\times T\Rr \\
                            &  (v_1,v_2,t,\lambda) & \mapsto & (p_2(v_1+tv_2),(t,t\cdot\exp{\lambda}))
\end{array}
\]
The isomorphism $di_{\D}^{-1}:\mathscr{A}\grVF\rightarrow \mathscr{A}\D_\varphi$ induced by $i_{\D}^{-1}$ over $\Omega_\mathscr{A}=\bigcup\Omega_i^j=N$ then implies the equality $p_{|N}=p_\varphi\circ di_{\D}^{-1}$ and $\A\grVF=\A\mathcal{F}_{M|N}$. 

Thus $\grVF$ is a quasi-graphoid integrating the Lie algebroid $\A\F$ over M. $\A\F$ being extremal, $\grVF$ is indeed the holonomy groupoid $\HolM$ of the foliation $(M,\F)$.
\end{demo}

Finally we get as a corollary a similar interpretation for the groupoid of a manifold with fibred corners:

\begin{corollary}
Let X be a manifold with fibred corners defined by a fibred decoupage $\eq_\phi=(M,(V_i,\phi_i)_{i\in I})$. Let $\F$ be the module of vector fields tangent to the fibrations of each face. Then 
$\gr \eq_\phi$ is the holonomy groupoid $\HolM$ of $(M,\F)$.
\end{corollary}

\begin{demo}
Let $\HolM$ be the holonomy groupoid of $(M,\F)$.
From previous proposition \ref{gr_holonomy}, for $i,j\in I$ the minimality of $\gr V_{i\phi}$ and $\gr V_{j\phi}$ implies the existence of morphisms $p_i$ and $p_j$ such that the following diagram commutes:
\[
\xymatrix{
\HolM \ar[d]_{p_j}\ar[r]^{p_i}           &  \gr V_{i\phi}
 \ar[d]^{s_i\oplus r_i} \\
\gr V_{j\phi} \ar[r]^{s_j\oplus r_j}    & M^2
}
\]

From the universal property of the fibered product $\gr \eq_\phi$ there thus exists a morphism of groupoids $u:\HolM\rightarrow \gr \eq_\phi$.
The minimality of $\HolM$ then implies $u$ to be an isomorphism and $\gr \eq_\phi\simeq \HolM$ is the holonomy groupoid of $(M,\F)$.
\end{demo}

The groupoid $\Gamma_\phi(X)$ of a manifold with fibred corners has therefore a natural geometric meaning as a singular foliation holonomy groupoid, its construction is independent of any particular description. It is an explicit example of a singular leaf space in the sense of \cite{ASk}.
This result allows the conceptual interpretation of {$\phi$-calculus} as the pseudodifferential calculus associated with the holonomy groupoid of the singular foliation defined by the manifold with fibred corners.


\bibliographystyle{amsalpha}
\bibliography{bibliographie}

\providecommand{\bysame}{\leavevmode\hbox to3em{\hrulefill}\thinspace}
\providecommand{\MR}{\relax\ifhmode\unskip\space\fi MR }
\providecommand{\MRhref}[2]{%
  \href{http://www.ams.org/mathscinet-getitem?mr=#1}{#2}
}
\providecommand{\href}[2]{#2}
\begin{thebibliography}{ALMP09}

\bibitem[ALMP09]{ALMP}
Pierre Albin, Eric Leichtnam, Rafe Mazzeo, and Paolo Piazza, \emph{The
  signature package on witt spaces, i. index classes}, arxiv:math.DG/09061568
  (2009).

\bibitem[AS63]{AS0}
M.~F. Atiyah and I.~M. Singer, \emph{The index of elliptic operators on compact
  manifolds}, Bull. Amer. Math. Soc. \textbf{69} (1963), 422--433.

\bibitem[AS68a]{AS}
\bysame, \emph{The index of elliptic operators. {I}}, Ann. of Math. (2)
  \textbf{87} (1968), 484--530.

\bibitem[AS68b]{AS3}
\bysame, \emph{The index of elliptic operators. {III}}, Ann. of Math. (2)
  \textbf{87} (1968), 546--604.

\bibitem[AS06]{ASk}
Iakovos Androulidakis and Georges Skandalis, \emph{The holonomy groupoid of a
  singular foliation}, arxiv:math.DG/0612370 (2006).

\bibitem[AS09]{ASk2}
I.~{Androulidakis} and G.~{Skandalis}, \emph{{Pseudodifferential calculus on a
  singular foliation}}, ArXiv e-prints (2009).

\bibitem[BHS91]{BHS}
J.-P. Brasselet, G.~Hector, and M.~Saralegi, \emph{Th\'{e}or\`{e}me de de rham
  pour les vari\'{e}t\'{e}s stratifi\'{e}es}, Ann. Global Analy. Geom.
  \textbf{9} (1991), no.~3, 211--243.

\bibitem[BP85]{BP}
B~Bigonnet and Jean Pradines, \emph{Graphe d'un feuilletage singulier}, C.R.
  cad. Sci. Paris \textbf{300} (1985), no.~13, 439--442.

\bibitem[Con79]{Coinc}
Alain Connes, \emph{Sur la th\'eorie non commutative de l'int\'egration},
  Alg\`ebres d'op\'erateurs (S\'em., Les Plans-sur-Bex, 1978), Lecture Notes in
  Math., vol. 725, Springer, Berlin, 1979, pp.~19--143.

\bibitem[Con82]{Cosf}
\bysame, \emph{A survey of foliations and operator algebras}, Operator algebras
  and applications, Part I (Kingston, Ont., 1980), Proc. Sympos. Pure Math.,
  vol.~38, Amer. Math. Soc., Providence, R.I., 1982, pp.~521--628.

\bibitem[Con94]{Concg}
\bysame, \emph{Noncommutative geometry}, Academic Press Inc., San Diego, CA,
  1994.

\bibitem[CR07]{Ca3}
Paulo Carrillo-Rouse, \emph{Indices analytiques {\`a} support compact pour des
  groupoides de {L}ie}, Th{\`e}se de {D}octorat {\`a} l'{U}niversit{\'e} de
  {P}aris 7 (2007).

\bibitem[CS84]{CS}
Alain Connes and Georges Skandalis, \emph{The longitudinal index theorem for
  foliations}, Publ. Res. Inst. Math. Sci. \textbf{20} (1984), no.~6,
  1139--1183.

\bibitem[Deb01]{Deb}
Claire Debord, \emph{Holonomy groupoids of singular foliations}, J.
  Differential Geom. \textbf{58} (2001), no.~3, 467--500.

\bibitem[DL09]{DL}
Claire Debord and Jean-Marie Lescure, \emph{K-duality for stratified
  pseudomanifolds}, Geometry and Topology \textbf{13} (2009), 49--86.

\bibitem[DLN06]{DLN}
Claire Debord, Jean-Marie Lescure, and Victor Nistor, \emph{Groupoids and an
  index theorem for conical pseudomanifolds}, arxiv:math.OA/0609438 (2006).

\bibitem[DLR11]{DLR}
C.~{Debord}, J.-M. {Lescure}, and F.~{Rochon}, \emph{{Pseudodifferential
  operators on manifolds with fibred corners}}, ArXiv e-prints (2011).

\bibitem[Gui12]{Gui}
Laurent Guillaume, \emph{G{\'e}om{\'e}trie non-commutative et calcul
  pseudodiff{\'e}rentiel sur les vari{\'e}t{\'e}s {\`a} coins fibr{\'e}s},
  Ph.D. thesis, {U}niversit{\'e} {P}aul {S}abatier {T}oulouse 3, 2012.

\bibitem[Hae84]{Haef}
Andr{\'e} Haefliger, \emph{Groupo\"\i des d'holonomie et classifiants},
  Ast\'erisque (1984), no.~116, 70--97, Transversal structure of foliations
  (Toulouse, 1982).

\bibitem[HS87]{HS}
Michel Hilsum and Georges Skandalis, \emph{Morphismes {$K$}-orient\'es
  d'espaces de feuilles et fonctorialit\'e en th\'eorie de {K}asparov
  (d'apr\`es une conjecture d'{A}. {C}onnes)}, Ann. Sci. \'Ecole Norm. Sup. (4)
  \textbf{20} (1987), no.~3, 325--390.

\bibitem[LMN00]{LMN}
Robert Lauter, Bertrand Monthubert, and Victor Nistor, \emph{Pseudodifferential
  analysis on continuous family groupoids}, Doc. Math. \textbf{5} (2000),
  625--655 (electronic).

\bibitem[LMN05]{LMN2}
Robert Lauter, Bertrand Monthubert, and Victor Nistor, \emph{Spectral
  invariance for certain algebras of pseudodifferential operators}, Journal of
  the Institute of Mathematics of Jussieu \textbf{4} (2005), no.~03, 405--442.

\bibitem[Mat70]{Mather}
J.N. Mather, \emph{Notes on topological stability}, Mimeographed Notes, 1970.

\bibitem[Mel93]{Mel}
Richard~B. Melrose, \emph{The {A}tiyah-{P}atodi-{S}inger index theorem},
  Research Notes in Mathematics, vol.~4, A K Peters Ltd., Wellesley, MA, 1993.

\bibitem[MM98]{MM}
Rafe Mazzeo and Richard~B. Melrose, \emph{Pseudodifferential operators on
  manifolds with fibred boundaries}, arxiv:math.DG/9812120 (1998).

\bibitem[Mon03]{Mont}
Bertrand Monthubert, \emph{Groupoids and pseudodifferential calculus on
  manifolds with corners}, J. Funct. Anal. \textbf{199} (2003), no.~1,
  243--286.

\bibitem[MP97]{MP}
Bertrand Monthubert and Fran{\c{c}}ois Pierrot, \emph{Indice analytique et
  groupo\"\i des de {L}ie}, C. R. Acad. Sci. Paris S\'er. I Math. \textbf{325}
  (1997), no.~2, 193--198.

\bibitem[MR05]{MRo}
R.~B. {Melrose} and F.~{Rochon}, \emph{{Index in K-theory for families of
  fibred cusp operators}}, ArXiv Mathematics e-prints (2005).

\bibitem[NWX99]{NWX}
Victor Nistor, Alan Weinstein, and Ping Xu, \emph{Pseudodifferential operators
  on differential groupoids}, Pacific J. Math. \textbf{189} (1999), no.~1,
  117--152.

\bibitem[Ren80]{Ren}
Jean Renault, \emph{A groupoid approach to {$C\sp{\ast} $}-algebras}, Lecture
  Notes in Mathematics, vol. 793, Springer, Berlin, 1980.

\bibitem[Tho69]{Thom}
Ren{\'e} Thom, \emph{Ensembles et morphismes stratifi{\'e}s}, Bull. Amer. Math.
  Soc. \textbf{75} (1969), 240--284.

\bibitem[Ver84]{Verona}
Andrei Verona, \emph{Stratified mappings-structure and triangulability},
  Lecture Notes in Mathematics \textbf{1102} (1984).

\bibitem[Win83]{Win}
H.~E. Winkelnkemper, \emph{The graph of a foliation}, Ann. Global Anal. Geom.
  \textbf{1} (1983), no.~3, 51--75.

\end{thebibliography}

\end{document}